\documentclass[12pt,a4paper,reqno]{amsart}
\usepackage{amsmath, amssymb, amscd, amsthm,latexsym,graphics,epsfig}

\title{Conical limit sets and continued fractions}

\author{Edward Crane}
\address{Department of Mathematics \\
University Walk \\
Bristol \\
BS8 1TW\\
United Kingdom}
\email{Edward.Crane@gmail.com}
\thanks{The first author was supported by a junior research fellowship at Merton College, Oxford, and by the University of Bristol.}

\author{Ian Short}
\address{Logic House \\
National University of Ireland, Maynooth \\
Maynooth\\
County Kildare\\
Ireland}
\email{Ian.Short@nuim.ie}
\thanks{The second author was supported by Science Foundation Ireland grant 05/RFP/MAT0003.}

\date{\today}


\newtheorem{theorem}{Theorem}
\newtheorem{lemma}[theorem]{Lemma}
\newtheorem{corollary}[theorem]{Corollary}
\newtheorem*{theorem*}{Theorem}

\newtheorem{proposition}[theorem]{Proposition}
\theoremstyle{definition}
\newtheorem{definition}[theorem]{Definition}
\newtheorem{problem}{Problem}

\numberwithin{theorem}{section}
\numberwithin{equation}{section}

\renewcommand{\leq}{\leqslant}
\renewcommand{\geq}{\geqslant}

\begin{document}

\newcommand{\Chat}{{\mathbb{C}_\infty}}
\newcommand{\Rhat}{{\mathbb{R}_\infty}}
\newcommand{\jj}{\ensuremath{\jmath}}

\renewcommand{\leq}{\leqslant}
\renewcommand{\geq}{\geqslant}
\renewcommand{\le}{\leqslant}
\renewcommand{\ge}{\geqslant}

\begin{abstract}

Inspired by questions of convergence in continued fraction theory,
Erd\H{o}s,  Piranian and Thron studied the possible sets of divergence for arbitrary sequences of M\"obius maps acting on the Riemann sphere, $S^2$. By identifying $S^2$ with the boundary of three-dimensional hyperbolic space, $H^3$, we show that these sets of divergence are precisely the sets that arise as conical limit sets of subsets of $H^3$. Using hyperbolic geometry, we give simple geometric proofs of the theorems of Erd\H{o}s, Piranian and Thron that generalise  to arbitrary dimensions. New results are also obtained about the class of conical limit sets, for example, that it is closed under locally quasisymmetric homeomorphisms. Applications are given to continued fractions. 
\end{abstract}

\maketitle

\section{Introduction}\label{S: intro}{\quad}\vspace{6pt}

In \cite{PT} George Piranian and Wolfgang Thron defined a class of subsets of the Riemann sphere $\mathbb{C}_\infty$, which they called \emph{sets of divergence}. A set $A \subseteq\mathbb{C}_\infty$ is a set of divergence if and only if there exists a sequence $\left(G_n\right)$ of M\"obius transformations such that the sequence $\left(G_n(z)\right)$ diverges for each $z\in A$ and converges for each $z \not\in A$. Paul Erd\H{o}s and Piranian continued the study of these sets in \cite{EP}, giving a geometric characterisation of the countable sets of divergence. The motivation for these papers was the theory of real and complex continued fractions. Associated to each complex continued fraction is a sequence of M\"obius transformations, and convergence of the continued fraction is equivalent to convergence of this sequence at $z=0$. For each $m \geq 1$, let $H^{m+1}$ denote $(m+1)$-dimensional hyperbolic space and $S^{m}$ denote its ideal boundary sphere. M\"{o}bius maps acting on $S^m$ extend uniquely to isometries of hyperbolic space $H^{m+1}$. Alan Beardon and others have exploited this isometric action to give simple geometric proofs of various classical results about convergence of complex continued fractions, avoiding explicit estimates on coefficients in matrix products. These proofs still apply in higher dimensions.

Suppose that a sequence $\left(z_n\right)$ in $H^{m+1}$ converges ideally to a limit point $x$ in $S^{m}$. The sequence is said to converge \emph{conically} if there is a geodesic ray $\gamma$ in $H^{m+1}$, landing at $x$, such that the points $z_n$ are all within a bounded hyperbolic distance from $\gamma$. For any subset $E \subseteq H^{m+1}$, the \emph{conical limit set} of $E$, denoted $\Lambda_c(E)$, consists of those points $x$ in $S^m$ for which there is a sequence of points in $E$ that converges conically to $x$.

\begin{proposition}\label{P: main connection}
A subset of $S^{m}$ is a set of divergence for some sequence of orientation preserving M\"{o}bius maps acting on $S^{m}$ if and only if it is the conical limit set of some subset of ${H}^{m+1}$.
\end{proposition}

Proposition~\ref{P: main connection} enables us to give purely geometric proofs of the known results about sets of divergence, and also some new results about these sets, without reference to explicit sequences of M\"{o}bius maps. Most of the proofs of \cite{EP,PT} involve explicit sequences of M\"obius transformations, hence these proofs are valid only in two-dimensions. In contrast, our geometric proofs are valid in all dimensions. Our first main result consists of generalisations of the topological results of \cite{EP,PT} to arbitrary dimensions. Let $\textup{CL}(m)$ be the class of all subsets of $S^{m}$ that arise as conical limit sets of subsets of $H^{m+1}$.
\begin{theorem}\label{T: main 1}
The following hold for each $m \geq 1$.
\begin{enumerate}
\item $\textup{CL}(m)$ is contained in, but not equal to the class of $G_{\delta\sigma}$ subsets of $S^m$.
\item The class of $G_\delta$ subsets of $S^m$ is contained in, but not equal to $\textup{CL}(m)$.
\item Let $Y$ be a sphere of codimension one in $S^m$. Then $\textup{CL}(m)$ contains all $G_{\delta\sigma}$ subsets of $Y$.
\item $\textup{CL}(m)$ is not invariant under homeomorphisms of $S^m$. 
\end{enumerate}
\end{theorem}

The geometrical point of view suggested a number of new results about conical limit sets. They are simple to prove, given the geometrical machinery available, but if the proofs were written in terms of sets of divergence of M\"{o}bius sequences they would seem less transparent.

\begin{theorem}\label{T: main 2}{\quad}
\begin{enumerate}
\item If $\left\{U_i\,: i \in I\right\}$ is a cover of $X$ by open subsets of $S^m$ then $X \in \textup{CL}(m)$ if and only if $X \cap U_i \in \textup{CL}(m)$ for each $i \in I$.
\item Two subsets $E$ and $F$ of  $S^m$ are respectively the conical limit set and the limit set of some subset of $\mathbb{H}^{m+1}$ if and only if $E \in \textup{CL}(m)$, $F$ is closed, $E$ is contained in $F$, and the interior of $F$ is contained in the closure of $E$.
\item $\textup{CL}(m)$ is preserved under quasiconformal homeomorphisms of $S^m$ for $m \ge 2$, and under quasisymmetric homeomorphisms in the case $m=1$. 
\end{enumerate}
\end{theorem}

Our final theorem is the generalisation to arbitrary dimensions of Erd\H{o}s and Piranian's characterisation of the countable sets of divergence of M\"{o}bius sequences acting on $\mathbb{C}_\infty$. The new feature of our work is  the geometrical interpretation of their theorem and proof; the generalisation to arbitrary dimension is a consequence. For any subset $E \subseteq S^{m}$, let $\overline{\textup{co}}(E)$ denote the hyperbolic convex hull of $E$ in $\mathbb{H}^{m+1}$. Then we define inductively a set $E^\chi$ for each ordinal $\chi$, as follows:
\begin{eqnarray*} 
 E^1 & = & E,  \\
 E^{\chi + 1} & =&  E^\chi \setminus \Lambda_c(\mathbb{H}^{m+1} \setminus \overline{\textup{co}}(E^\chi)), \quad\text{and}\\
 E^{\chi} & = & \bigcap_{\psi < \chi} E^\psi, \quad \text{for a limit ordinal $\chi$}.
\end{eqnarray*}  
\begin{theorem}\label{T: main 3}
For $m \geq 1$, a countable set $E \subseteq S^m$ is in $\textup{CL}(m)$ if and only if there exists an ordinal $\chi$ such that $E^{\chi} = \emptyset$. 
\end{theorem}

The authors are grateful for the detailed comments of the referee, which have been used to improve the paper immensely.

\section{Background and Definitions}\label{S: background}

\subsection{M\"{o}bius maps and hyperbolic space}\label{SS: hyperbolic space}{\quad}\vspace{6pt}

 Much of this paper will be concerned with groups of M\"{o}bius maps and with the geometry of $H^m$ ($m\geq 2$), the $m$-dimensional hyperbolic space. In this section we set out some basic facts about M\"{o}bius maps and models of $H^m$, in order to establish notation. We will use standard properties of hyperbolic geometry without comment, and for more detail the reader may wish to consult \cite{BeardonDG} or \cite{Rat}.

The hyperbolic space $H^m$ has a Riemannian metric $ds^2$ with distance function $\rho(x,y)$. The space $H^m$ is the unique connected and simply-connected $m$-dimensional Riemannian manifold of constant curvature $-1$.  The full group of isometries of ${H}^m$ is denoted $\textup{Isom}\left({H}^m\right)$. The subgroup of orientation-preserving isometries is denoted $\textup{Isom}^+\left({H}^m\right)$. 


  We will use two standard conformal models of the metric space $\left({H}^m, \rho\right)$: the ball model and the upper half-space model. In both models we will denote the hyperbolic distance by $\rho$, but we will always make it clear which model we are using when giving explicit constructions. Before introducing these two models of $H^m$, it is necessary to define M\"obius transformations in several dimensions.
 
 The space $\mathbb{R}^m_\infty$ is the one-point compactification of Euclidean space $\mathbb{R}^m$; it consists of $\mathbb{R}^m$ together with a single extra point $\infty$. The M\"{o}bius group $\textup{M\"ob}(m)$ is the full group of conformal homeomorphisms of $\mathbb{R}^m_\infty$. It is generated by reflections in planes and spherical inversions.


 The Poincar\'e ball model consists of the open unit ball 
\[
\mathbb{B}^m = \{ x \in \mathbb{R}^m \,:\,|x|<1\}\,,
\] 
equipped with the hyperbolic metric given by the Riemannian metric 
\[ 
ds^2 = \frac{4\,|dx|^2}{\left(1-|x|^2\right)^2}\,.
\] 
The geodesics in the ball model are arcs of lines and circles in $\mathbb{R}^m$ that cut the boundary sphere orthogonally.  The group $\textup{Isom}\left(\mathbb{B}^m\right)$ of hyperbolic isometries of $\mathbb{B}^m$ is the subgroup of $\textup{M\"{o}b}(m)$ consisting of maps that preserve $\mathbb{B}^m$. This group is generated by reflections that fix $0$ and spherical inversions in hyperspheres that meet the unit sphere orthogonally. 

 Note that for $m < n$ the inclusion of $\mathbb{R}^m$ as $\mathbb{R}^m \times \{0\} \subseteq\mathbb{R}^m \times \mathbb{R}^{n-m} = \mathbb{R}^n$ induces an isometric embedding of hyperbolic metrics, $\mathbb{B}^m \hookrightarrow \mathbb{B}^n$. Any isometric embedding of ${H}^m$ into ${H}^n$ can be mapped to this standard one by a suitable identification of ${H}^n$ with $\mathbb{B}^n$. Each isometry of ${H}^m$ can be extended to an isometry of ${H}^n$, though not in a unique way. 

The \emph{ideal boundary} of ${H}^m$ is a sphere $S^{m-1}$, which we sometimes denote $\partial H^m$. It can be defined intrinsically in terms of equivalence classes of geodesic rays, but it is easiest to explain using the ball model, where it is represented  by the unit sphere $\mathbb{S}^{m-1}$ which is the boundary in $\mathbb{R}^m$ of $\mathbb{B}^m$. A sequence of points  $(z_n)$ in ${H}^m$ \emph{converges ideally} to a point $x$ in $\partial H^m$ if the representative points in $\mathbb{B}^m$ converge in the Euclidean topology to a point of $\mathbb{S}^{m-1}$; we write this simply as $z_n \to x$ as $n \to \infty$, even though $(z_n)$ does not converge in the hyperbolic metric. The isometries of $H^m$ extend to homeomorphisms of $H^m \cup \partial H^m$. 


We use the notation $\jj$ to a refer to a distinguished point in $H^m$; in the ball model of hyperbolic space we always choose $\jj$ to be the origin in $\mathbb{R}^m$. Let $\gamma$ be a hyperbolic geodesic in $\mathbb{B}^m$ with endpoints $x,y \in \mathbb{S}^{m-1}$, and let $\rho(\jj,\gamma)$ denote the hyperbolic distance of $\jj$ from $\gamma$. The next equation provides a useful link between the hyperbolic metric in $\mathbb{B}^m$ and the restriction of Euclidean distance in $\mathbb{R}^m$ to the unit sphere $\mathbb{S}^{m-1}$:
\begin{equation}\label{E: key geometry}
 |x-y|\,=\,2/\cosh \rho(\jj,\gamma)\,.
\end{equation} 
In $\mathbb{B}^m$, a \emph{horoball} based at $x \in \mathbb{S}^{m-1}$ is an open Euclidean ball of radius less than $1$ which is internally tangent to $\mathbb{S}^{m-1}$ at $x$. The boundary of a horoball is called a \emph{horosphere}. 
  
Next, we describe the upper half-space model of ${H}^m$, which is
\[
 \mathbb{H}^m \, = \, \{ \left(x_0, \dots, x_{m-1}\right) \in \mathbb{R}^m \, : \, x_0 > 0\}\,,
\]
 equipped with the metric \[ds^2 = \frac{|dx|^2}{x_0^2}\,.\] 
There is a standard isometric mapping $\mathbb{B}^m \to \mathbb{H}^m$ given by a spherical inversion in the point $(-1, 0, \dots, 0)$ in $\mathbb{R}^m_\infty$; we will not need to use it explicitly. Our distinguished point $\jj \in H^m$ is represented in $\mathbb{H}^m$ by the point $e_0 = (1,0,\dots,0)$.
 
 The geodesics in $\mathbb{H}^m$ are arcs of lines and circles in $\mathbb{R}^m$ that cut the boundary plane orthogonally.  The ideal boundary $\partial {H}^m$ is represented by $\mathbb{R}^{m-1}_\infty$, which is the boundary of $\mathbb{H}^m$ within $\mathbb{R}^m$ together with the point at infinity. A sequence $\left(z_n\right)$ in $\mathbb{H}^m$ converges ideally to a point $x$ in $\mathbb{R}^{m-1}_\infty$ if and only if $\left(z_n\right)$ converges to $x$ in the topology of the one-point compactification $\mathbb{R}^m_\infty$. 

 In the special case of ${H}^3$, using the upper half-plane model we can identify the ideal boundary with the Riemann sphere $\mathbb{C}_\infty = \mathbb{C} \cup \{\infty\}$, and the orientation-preserving M\"{o}bius maps act on $\mathbb{C}_\infty$ as fractional linear transformations \[z \mapsto \frac{az+b}{cz+d}\,, \quad ad-bc \neq 0\,.\] Thus $\textup{Isom}^+\left({H}^3\right) \cong \textup{PSL}(2,\mathbb{C})$. The subgroup $\textup{PSL}(2,\mathbb{R})$ acts as the group of orientation-preserving isometries of the standard hyperbolic plane $\mathbb{H}^2$ whose ideal boundary is $\mathbb{R}_\infty \subseteq\mathbb{C}_\infty$.

\subsection{Conical limits}\label{SS: hyperbolic space conical limit points}{\quad}\vspace{6pt}


We denote the hyperbolic distance between a point $z$ and a set $E$ in ${H}^m$ by \[\rho(z,E) = \inf_{w \in E} \rho(z,w)\,.\] For $\alpha > 0$, the $\alpha$-neighbourhood of $E$ is defined by \[N_\alpha(E) = \{z \in {H}^m\,:\, \rho(z,E) < \alpha\}\,.\] We denote by $[x,z]$ the geodesic segment from $x$ to $z$, where $x$ and $z$ may be points of $H^m$ or of its ideal boundary $S^{m-1} = \partial H^m$; any ideal endpoints are not included in $[x,z]$, so that $[x,z]$ is a subset of $H^m$.

\begin{definition}
 Let $\left(z_n\right)$ be a sequence of points in ${H}^m$. For any $\alpha > 0$ we say that the ideal boundary point $x \in S^{m-1}$ is the \emph{$\alpha$-conical limit point of $\left(z_n\right)$} when
 \begin{enumerate}
 \item $z_n \to x$ as $n \to \infty$, and
 \item $z_n \in N_{\alpha}([\jj,x])$ for all sufficiently large $n$.
 \end{enumerate}
 We say that $(z_n)$ \emph{converges conically} to $x$ if there is $\alpha>0$ such that $x$ is the $\alpha$-conical limit point of $(z_n)$. The boundary point $x$ is an $\alpha$-conical limit point of a set $E \subseteq H^m$ if and only if $x$ is the $\alpha$-conical limit of some sequence of points $(z_n)$ in  $E$. The set of $\alpha$-conical limit points of $E$ is denoted $\Lambda_c^{\alpha}(E)$. Finally, $x$ is a conical limit point of $E$ if it is an $\alpha$-conical limit point of $E$ for some positive $\alpha$, and the conical limit set of $E$ is the set of all conical limit points of $E$: 
\begin{equation}\label{E: conical limit set formula}
\Lambda_c(E) = \bigcup_{\alpha>0} \Lambda_c^{\alpha}(E) \,.
\end{equation}
\end{definition}
\noindent If $E$ is a bounded subset of hyperbolic space, then $\Lambda_c(E) = \emptyset$.

\begin{figure}[ht]
\centering
\includegraphics[scale=0.7]{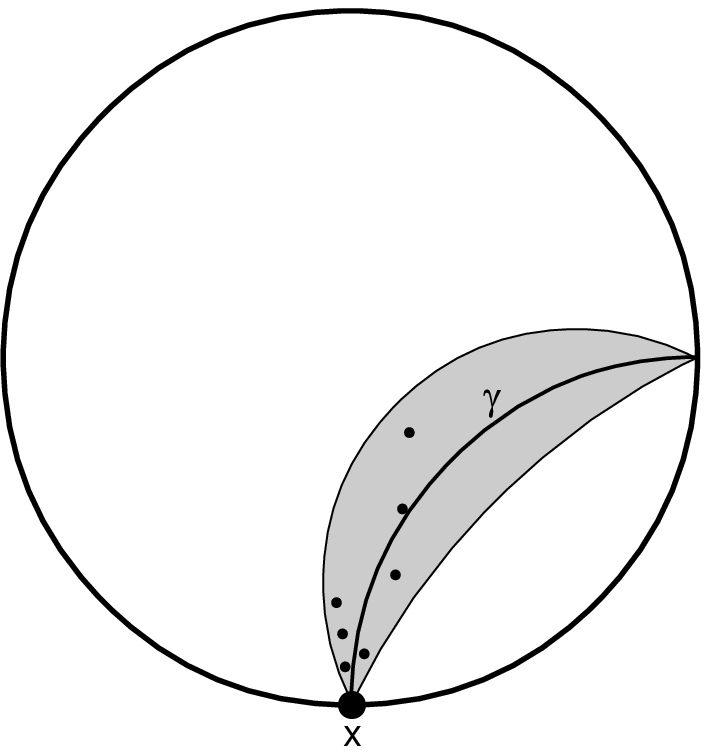}
\caption{}
\label{F: conicalConvergence}
\end{figure}

\begin{lemma}\label{L: any geodesic will do}
Suppose that $(z_n)$ converges $\alpha$-conically to $x$. Then for each $\alpha' > \alpha$ and each geodesic $\gamma$ ending at $x$, we have $z_n \in N_{\alpha'}(\gamma)$ for $n$ sufficiently large. 
\end{lemma}
\begin{proof}
This is a simple consequence of the fact that if two geodesics $\gamma$ and $\gamma'$ both end at $x$, then  $\rho(z, \gamma') \to 0$ as $z \to x$ along $\gamma$.
\end{proof}
In Figure~\ref{F: conicalConvergence} a sequence is shown converging conically to a limit $x$.
\begin{corollary}
 If $g \in \textup{Isom}\left({H}^m\right)$ then 
\[
\Lambda_c(g(E))\, = \,g\left(\Lambda_c(E)\right)\,.
\]
\end{corollary}
\begin{definition} For $m \ge 2$ we define the class of conical limit sets in the sphere $S^{m-1} = \partial H^m$ as follows:
 \[\textup{CL}(m-1) \,=\, \left\{ A \subseteq \partial H^m \,:\, A = \Lambda_c(E) \; \text{ for some } \; E \subseteq H^m \right \}. \]
We use the same notation to refer to the class of conical limit sets in the ideal boundary in a specific model of hyperbolic space, for example, the class of subsets of $\mathbb{R}^{m-1}_\infty$ that arise as conical limit sets of subsets of $\mathbb{H}^m$.
\end{definition}
 
Next, we give an alternative description of $\Lambda_c(E)$, which will be useful for investigating the structure of conical limit sets. 

\begin{definition}
 Let $w \in H^m$ and $v \in H^m \cup \partial H^m$. The $\alpha$-shadow of $w$ from $v$ is
\[\textup{Shad}_v^\alpha(w) = \{ x \in \partial H^m \, : \, \rho(w, [v,x]) < \alpha \} = \{ x \in \partial H^m \, : \, w \in N_\alpha([v,x]) \}\,. \]
\end{definition}
\noindent Imagining that light travels along hyperbolic geodesics, we may think of $\textup{Shad}_v^{\alpha}(w)$ as the shadow cast on the ideal boundary by the open ball of hyperbolic radius $\alpha$ centred on $w$, when illuminated from the point $v$. The shadow is an open ball in the spherical metric of $\partial H^m$.

\begin{lemma}\label{L: shadow descriptions of conical limit sets}\text{}
Let $E$ be a subset of $H^m$.
\begin{enumerate}
\item Let $\overline{B_\rho(\jj,r)}$ denote the closed ball or hyperbolic radius $r$ about $\jj$. Then
 \[
\Lambda_c^\alpha(E) = \bigcap_{r > \alpha \sqrt{3}}\; \bigcup_{w \in E \setminus \overline{B_{\rho}(\jj,r)}} \textup{Shad}_\jj^\alpha(w)\,.
\]
The shadows in this formula have diameter less than $\pi/3$ in the spherical metric.
\item For $x \in \partial H^m$, let $\textup{HB}(x)$ denote the family of horoballs based at $x$. For each $\epsilon \in (0,\alpha)$,
\[
\Lambda_c^{\alpha - \epsilon}(E) \setminus \{x \} \,\subseteq\, \bigcap_{D \in \textup{HB}(x)}\, \bigcup_{w \in E \setminus D} \textup{Shad}_x^\alpha(w)\, \subseteq \, \Lambda_c^{\alpha + \epsilon}(E) \setminus \{x\}\,.\]
\end{enumerate}
\end{lemma}
\begin{proof}
The first statement is a simple consequence of the definitions. The second statement follows from Lemma~\ref{L: any geodesic will do}.
\end{proof}

We frequently apply Lemma~\ref{L: shadow descriptions of conical limit sets}  (ii) in the upper half-space model of hyperbolic space, with $x=\infty$, because both the neighbourhoods $N_\alpha(\gamma)$ and the $\alpha$-shadows have simple geometric interpretations in this model. Let $\delta_x$ be the `vertical' geodesic in  $\mathbb{H}^{m}$  with endpoints $x\in\mathbb{R}^{m-1}$ and $\infty$. The next lemma says that the $\alpha$-neighbourhood $N_\alpha(\delta_x)$ is an infinite cone with axis $\delta_x$ and vertex at $x$. This makes the model $\mathbb{H}^{m}$ particularly useful for examining conical limit sets, since it allows us to describe the conical limit set of a subset of $\mathbb{H}^{m}$ entirely in terms of Euclidean cones. 

\begin{lemma}\label{L: cones}\cite[\S 7.20]{BeardonDG}
In the upper half-space model, the $\alpha$-neighbourhood $N_{\alpha}(\delta_x)$ is the interior of a cone, rotationally symmetric about the geodesic $\delta_x$, with vertex at $x$. The angle at $x$ between $\delta_x$ and the boundary of the cone is $\theta$, where $\cos \theta \, \cosh \alpha = 1$. 
\end{lemma}

\begin{figure}[ht]
\centering
\includegraphics[scale=0.7]{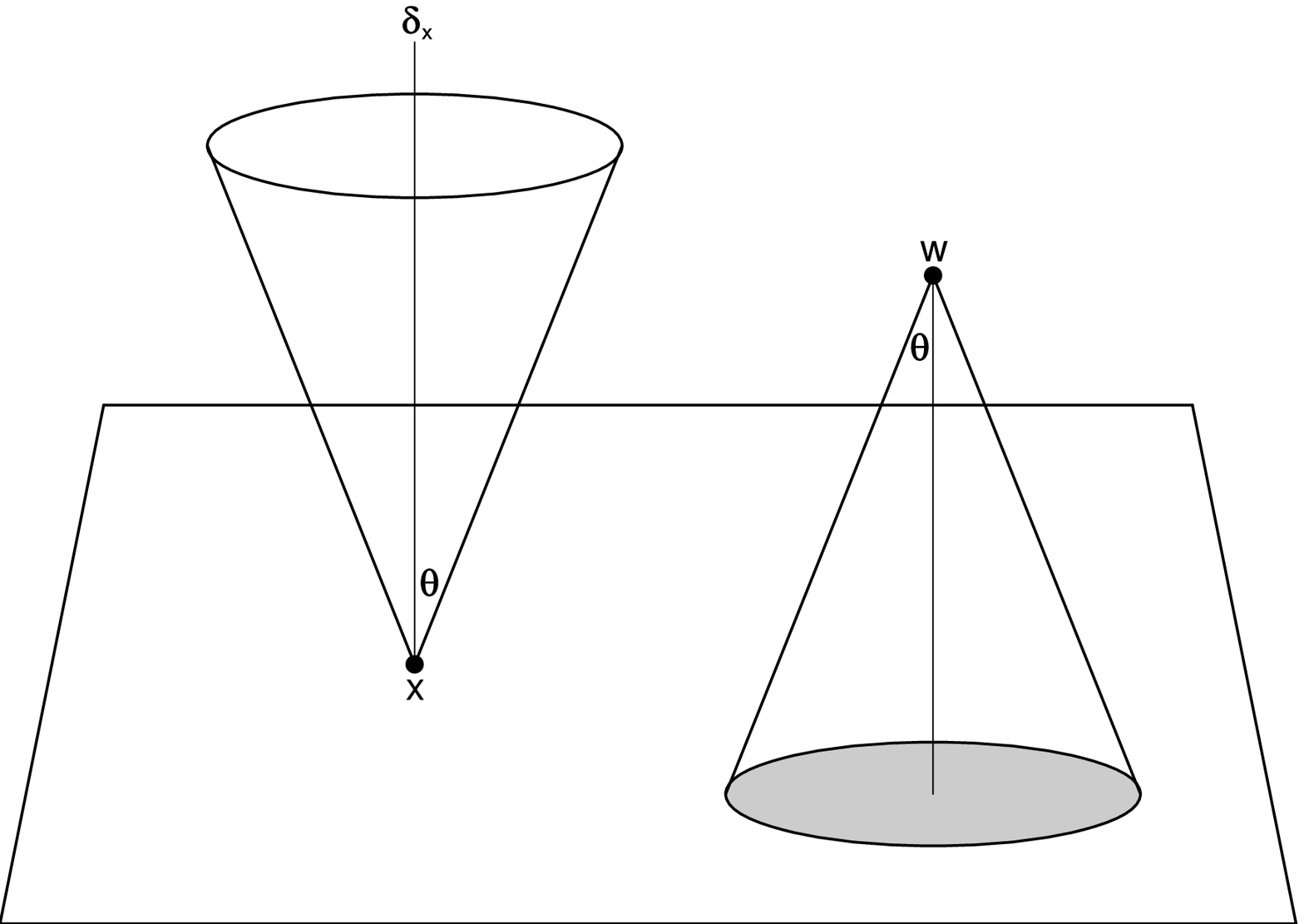}
\caption{}
\label{F: cone}
\end{figure}


\begin{lemma}\label{L: alpha shadows from infinity}
 Let $w = (t,v_1, \dots, v_{m-1}) \in \mathbb{H}^{m}$, where $v = (v_1, \dots, v_{m-1}) \in \mathbb{R}^{m-1}$ and $t>0$. Then we have
  \[
\textup{Shad}_\infty^\alpha(w) = \{ x \in \mathbb{R}^{m-1} \, : \, |x - v| < t \sinh \alpha\,\}\,.
\]
\end{lemma}

The left hand cone of Figure~\ref{F: cone} depicts the neighbourhood $N_\alpha(\delta_x)$.  Lemma~\ref{L: cones} says that $\rho(w, \delta_u) < \alpha$ if and only if $u \in \textup{Shad}_\infty^\alpha(w)$. The shaded region of the right hand cone of Figure~\ref{F: cone} depicts the $\alpha$-shadow from $\infty$ of the point $w$.

The Hausdorff distance $\rho_H(E,F)$ between two subsets $E$ and $F$ of $H^m$ is defined by the equation 
\[ 
\rho_H(E,F) = \max\left( \sup_{z \in E} \rho(z,F),\, \sup_{w \in F} \rho(w,E) \right)\,.
\]
The Hausdorff distance is not a metric, although its restriction to the set of compact subsets of ${H}^m$ is a metric.
\begin{lemma}\label{L: finite Hausdorff distance}
 Let $E, F \subseteq {H}^m$. If $\rho_H(E,F) < \infty$ then $\Lambda_c(E) = \Lambda_c(F)$.
\end{lemma}

\begin{proof}
Let $M=\rho_H(E,F)$, and suppose $x \in \Lambda_c(E)$. Then there exists a sequence $\left(z_n\right)$ of points in $E$ that converges to $x$ inside $N_\alpha(\gamma)$, for some geodesic $\gamma$ ending at $x$ and some $\alpha>0$. For each $n$, we can pick $w_n \in F$ such that $\rho\left(z_n, w_n\right) \leq M+1$, so the sequence $\left(w_n\right)$ converges to $x$, and each $w_n$ lies in $N_{\alpha + M+1}(\gamma)$. Therefore $x$ is a conical limit point of $F$. We have shown that $\Lambda_c(E) \subseteq \Lambda_c(F)$ and the reverse inclusion follows by exchanging $E$ and $F$ in the argument.
\end{proof}

A well-known result which is an immediate application of Lemma~\ref{L: finite Hausdorff distance} is that the conical limit set of an orbit of a Kleinian group $\Gamma$ does not depend on the choice of orbit. The conical limit set plays an important role in the theory of Kleinian groups, being a natural $\Gamma$-invariant subset of the limit set of $\Gamma$.

The next few results concern the conical limit sets of countable subsets of $H^m$. We will use the abbreviation $\Lambda_c\left(w_n\right)$ for $\Lambda_c\left(\left\{w_n: n \in \mathbb{N}\right\}\right)$.

\begin{definition}\label{D: standard and escaping}
A sequence $\left(w_n\right)$ of points in ${H}^m$ is an \emph{escaping} sequence if $\rho(w_n, \jj) \to \infty$ as $n \to \infty$. Equivalently, an escaping sequence is a sequence in $H^m$ that  eventually leaves every compact subset of $H^m$. A sequence $\left(w_n\right)$ of points in the model $\mathbb{H}^m$ is a \emph{standard} sequence if the $x_0$-co-ordinate of $w_n$ converges to $0$ as $n \to \infty$. Equivalently, a standard sequence is a sequence in $H^m$ that eventually leaves every horoball based at $\infty$.
\end{definition}
In particular, every standard sequence is an escaping sequence. 

\begin{lemma}\label{L: escaping and standard sequences}{\quad}
\begin{enumerate}
\item Let $E$ be any subset of ${H}^m$. Then there exists an escaping sequence $\left(w_n\right)$ in ${H}^m$ such that $\Lambda_c(E) =  \Lambda_c\left(w_n\right)$.
\item Let $E$ be any subset of $\mathbb{H}^m$. Then there exists a standard sequence $\left(w_n\right)$ in $\mathbb{H}^m$ such that $\mathbb{R}^{m-1} \cap \Lambda_c(E) =  \Lambda_c\left(w_n\right)$.
\end{enumerate} 
\end{lemma}
\begin{proof}
First we prove (i). If $\Lambda_c(E) = \emptyset$, take $\left(w_n\right)$ to be an orbit of the group generated by a parabolic M\"{o}bius map $M$; this converges to the fixed point of $M$ but no subsequence converges conically. Now suppose $\Lambda_c(E) \neq \emptyset$. Take any tessellation of ${H}^m$ by a countable sequence of sets $(X_n)$ with uniformly bounded diameters. For example, we could use the tessellation by translates of a fundamental polyhedron for some co-compact Kleinian group. Pass to the subsequence of those $X_n$ which contain a point of $E$, and for each of these take $w_n$ to be an arbitrary point in $X_n$. This subsequence is infinite because $\Lambda_c(E) \neq \emptyset$. The sequence $\left(w_n\right)$ is escaping because each compact set in ${H}^m$ meets only finitely many of the $X_n$. The Hausdorff distance between $E$ and $\{w_n \,:\,n \in \mathbb{N}\}$ is finite because the $X_n$ have uniformly bounded diameters. Thus Lemma~\ref{L: finite Hausdorff distance} shows that $\Lambda_c(E) =  \Lambda_c\left(w_n\right)$.

For part (ii), in the case where $\mathbb{R}^{m-1} \cap \Lambda_c(E) \neq \emptyset$, we first choose a sequence $\left(w_n\right)$ as above such that $\Lambda_c(E)=\Lambda_c(w_n)$. We express each $w_n$ as $(t_n, v_n)$ with $t_n >0$ and $v_n \in \mathbb{R}^{m-1}$, and pass to the subsequence of those $w_n$ for which $t_n < 1/(1 + |v_n|)$. This yields a standard sequence with the same conical limit set as $\left(w_n\right)$, except that the point $\infty$ is removed.
\end{proof}

Putting together Lemma~\ref{L: shadow descriptions of conical limit sets} and Lemma~\ref{L: escaping and standard sequences}, we can now express arbitrary conical limit sets in terms of countable unions and intersections of open balls:

\begin{corollary}\label{C: conical limit set formulae}
For any conical limit set $A \in \textup{CL}(m-1)$, there exists an escaping sequence $\left(w_n\right)$ in $H^m$ such that
\[
A = \bigcup_{\alpha=1}^{\infty} \bigcap_{p=1}^\infty\bigcup_{n=p}^\infty \textup{Shad}_\jj^{\alpha}(w_n)\,.
\]
For each subset $A \subseteq \mathbb{R}^{m-1}$ which is a conical limit set of a subset of $\mathbb{H}^m$, there exists a standard sequence $\left(w_n\right)$ in $\mathbb{H}^m$ such that
\[
A = \bigcup_{\alpha=1}^{\infty} \bigcap_{p=1}^\infty\bigcup_{n=p}^\infty \textup{Shad}_\infty^{\alpha}(w_n)\,.
\]
\end{corollary}

\section{Sets of divergence and conical limit sets}\label{S: Mobius}

\subsection{Convergence of sequences of M\"obius transformations}\label{SS: Mobius convergence}{\quad}\vspace{6pt}

The purpose of this section is to explain the connection between pointwise convergence of a sequence of M\"obius transformations $\left(G_n\right)$ on the interior of hyperbolic space to a limit $x$ in $S^{m-1}$, and convergence of $\left(G_n\right)$ to $x$ on the boundary of hyperbolic space.
\begin{definition}\label{D: general convergence}
 We say that a sequence of M\"obius transformations $\left(G_n\right)$ in $\textup{Isom}(H^m)$ \emph{converges generally} if the sequence $\left(G_n(\jj)\right)$ converges ideally to some point $x$ in the ideal boundary sphere $S^{m-1}$.
\end{definition}
To make this concrete in the ball model of ${H}^m$, we consider the elements $G_n$ as M\"{o}bius maps acting on the closed unit ball. The sequence $\left(G_n\right)$ converges generally to $x$ if and only if  $\left(G_n\right)$ converges pointwise and hence locally uniformly to $x$ on the open unit ball. This implies that the sequence $\left(G_n\right)$  leaves every compact subset of the topological group $\textup{Isom}\left({H}^m\right)$, so it has no convergent subsequence in the topology of the group. 

The notion of general convergence was introduced to continued fraction theory by Lisa Jacobsen (now Lisa Lorentzen) in \cite[p.480]{Jac1}. She expressed it in terms of the action of $G_n$ on the boundary sphere, as we explain in the following proposition. 
\begin{proposition}\label{P: two notions of general convergence}\textup{(Aebischer \cite[Theorem 4.3]{Aeb}, see also \cite[Theorem 6.6]{Be1}).}\\
A sequence of maps $(G_n)$ in $\textup{Isom}\left(\mathbb{B}^m\right)$ converges generally to a point $x$ in $\mathbb{S}^{m-1}$ if and only if there exist two sequences of points $\left(u_n\right)$ and $\left(v_n\right)$ in $\mathbb{S}^{m-1}$ such that
\[ \lim_{n \to \infty} G_n\left(u_n\right) = \lim_{n \to \infty} G_n(v_n) = x,\quad \inf \left(|u_n - v_n|\right) > 0 .\]
\end{proposition}

In particular, if there are two distinct points $u,v \in S^{m-1}$ such that $\left(G_n(u)\right)$ and $\left(G_n(v)\right)$ converge to the same limit, then the sequence $\left(G_n\right)$ converges generally. The converse is false: $\left(G_n\right)$ may converge generally yet diverge at every point of $S^{m-1}$. To see this, choose any sequence $(x_n)$ in $\mathbb{B}^{m+1}$ that converges ideally, say to $x \in \mathbb{S}^{m}$, and pick $H_n \in \textup{Isom}\left(\mathbb{B}^{m+1}\right)$ such that $H_n(\jj) = x_n$. The stabiliser $\textup{Stab}(x_n)$ acts transitively on $\mathbb{S}^{m}$, and its elements are Lipschitz with respect to the spherical metric (with a constant that only depends on $\rho(0,x_n)$). Therefore we can find a sequence $M_{n,1}, \dots, M_{n,k(n)}$ of elements of $\textup{Stab}(x_n)$ with the property that for each $y \in \mathbb{S}^{m}$, the set $\left\{M_{n,i}(y)\,: 1 \le i \le k(n)\right\}$ contains a point within spherical distance $1/n$ of each point of $\mathbb{S}^m$. Now consider the M\"{o}bius sequence $\left(G_n\right)$ given explicitly by
\[M_{1,1} \circ H_1, \dots , M_{1,k(1)} \circ H_1, M_{2,1} \circ H_2, \dots, M_{2,k(2)}\circ H_2, M_{3,1}\circ H_3, \dots, M_{3,k(3)}\circ H_3, \dots \,.\] 
 This sequence $\left(G_n\right)$ converges generally, but for each point $y$ in $\mathbb{S}^m$, the sequence $\left(G_n(y)\right)$ is dense in $\mathbb{S}^m$, so it certainly diverges.

\begin{proposition}\label{P: Aebischer}\textup{(Aebischer \cite[Theorem 5.2]{Aeb}).}\\
Let $(G_n)$ be a sequence of M\"{o}bius maps that converges generally to a point $x$ in $S^{m-1}$. Then for each $z \in S^{m-1}$, the sequence $\left(G_n(z)\right)$ converges to $x$ if and only if $z \not\in \Lambda_c\left(G_n^{-1}(\jj)\right)$. 
\end{proposition}
Proposition~\ref{P: Aebischer} is well known to Kleinian group theorists. Aebischer's proof used Euclidean geometry and a Euclidean definition of the notion of conical limits. Starting from our hyperbolic definition of conical limits we can give a  simple proof using the fact that $\jj$ is far from a hyperbolic geodesic if and only if the ideal endpoints of that geodesic are close together in the spherical metric on the ideal boundary sphere. 

\begin{figure}[ht]
\centering
\includegraphics[scale=1.0]{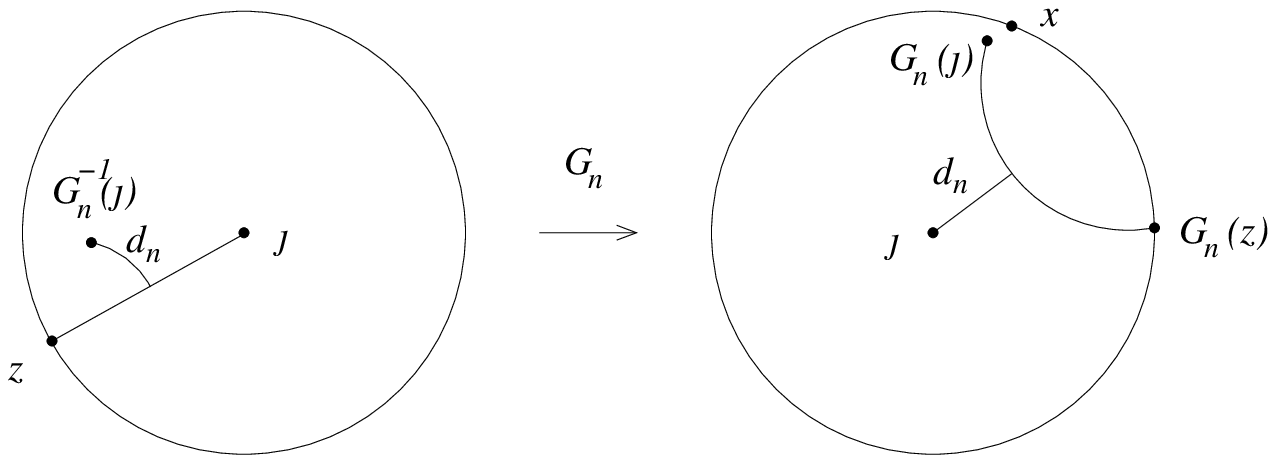}
\caption{}
\label{F: Aebischer}
\end{figure}

\begin{proof}[Proof of Proposition~\ref{P: Aebischer}]
Define $d_n = \rho\left(\jj, [G_n(\jj),G_n(z)]\right)$. Since $G_n(\jj) \to x$ as $n \to \infty$, equation~(\ref{E: key geometry}) shows that $G_n(z) \to x$ if and only if $d_n \to \infty$. Now we show that $d_n \not\to \infty$ if and only if $z \in \Lambda_c\left(G_n^{-1}(\jj)\right)$.  The hypothesis of general convergence implies that $\left(G_n^{-1}(\jj)\right)$ is an escaping sequence, because $\rho\left(G_n^{-1}(\jj),\jj\right) = \rho\left(\jj, G_n(\jj)\right) \to \infty$ as $n \to \infty$. Also, since $G_n$ is a hyperbolic isometry, we have $d_n = \rho\left(G_n^{-1}(\jj), [\jj,z]\right)$.
\end{proof}
  
\begin{lemma}\label{L: general convergence almost implies unique limit}
Suppose that $\left(G_n\right)$ converges generally to a point $x$. If for some element $u$ of  $S^{m-1}$ we have $G_n(u) \to y$, $y \neq x$, then the sequence $(G_n^{-1}(\jj))$ converges conically to $u$. In particular, $\left(G_n(u)\right)$ can converge to a point other than $x$ for at most one point $u$ in $S^{m-1}$, and if it does so then $G_n(z)\rightarrow x$ for each $z\neq u$. 
\end{lemma}
\begin{proof}
Since $\rho(\jj,[G_n(j),G_n(u)])\to \rho(\jj,[x,y])$ as $n\to\infty$ , the sequence with $n$th term $\rho(\jj,[G_n(\jj),G_n(u)])$ is bounded. But $\rho(G_n^{-1}(\jj),[\jj,u])=\rho(\jj,[G_n(\jj),G_n(u)])$ and $\left( G_n^{-1}(\jj)\right)$ is an escaping sequence, therefore $\left( G_n^{-1}(\jj)\right)$ converges conically to $u$.
 \end{proof}


\subsection{Proof of Proposition~\ref{P: main connection}}\label{SS: sets of divergence}{\quad}\vspace{6pt}

Piranian and Thron defined sets of divergence for sequences in $\textup{PSL}(2,\mathbb{C})$ acting as M\"{o}bius maps on the Riemann sphere. We now extend the definition in the obvious way to arbitrary dimensions. Let $m \ge 2$.
\begin{definition}
Let $\left(G_n\right)$ be any sequence in $\textup{Isom}^+\left(H^m\right)$, considered also to act on the ideal boundary sphere $\partial H^m = S^{m-1}$. The \emph{set of divergence} of $\left(G_n\right)$ is the subset of $S^{m-1}$ consisting of all points $z \in S^{m-1}$ such that the sequence $G_n(z)$ does not converge.
\end{definition}
Proposition~\ref{P: main connection} states that a subset $A$ of $S^{m-1}$ is a set of divergence for some sequence of orientation preserving M\"{o}bius maps acting on $S^{m-1}$ if and only if $A \in \textup{CL}(m-1)$. 
\begin{proof}
We may assume that $S^{m-1}\setminus A$ contains at least two points, since it is easy to check that the entire sphere and the complement of a single point arise both as sets of divergence and conical limit sets.

Let $(G_n)$ be a sequence of M\"obius maps in $\textup{Isom}\left({H}^m\right)$ that diverges only on the set $A$. We must show that $A$ is a conical limit set. There are two situations to consider. The first is that $(G_n)$ converges to a constant $x$ in $S^{m-1}$ on the complement of $A$. In this case $(G_n)$ converges generally to $x$ by Proposition~\ref{P: two notions of general convergence}. Then Proposition~\ref{P: Aebischer} shows that $A=\Lambda_c(G_n^{-1}(\jj))$. The second possibility is that $(G_n)$ does not converge to a constant on the complement of $A$. In this case it is known (\cite{Be2}) that $A$ is equal to either $\emptyset$, a singleton or the complement of a closed $k$-dimensional ball in ${S}^{m-1}$, for some $k<m-1$. It is straightforward to verify directly that these sets are all conical limit sets, although this will also follow from part (iii) of Theorem~\ref{T: main 1}.

For the converse we use the ball model of hyperbolic space. Let $A$ be a conical limit set containing at least two elements. Then there is a sequence $(z_n)$ in $\mathbb{B}^{m}$ with $|z_n|\to 1$ as $n\to\infty$, such that $\Lambda_c(z_n) = A$. Choose orientation preserving M\"obius transformations $\left(G_n\right)$ that satisfy $G_n^{-1}(\jj)=z_n$ and $G_n(\jj)=|z_n|e_0$. Then $\left(G_n\right)$ converges to $e_0$ on $\mathbb{S}^{m-1}\setminus A$ and $\left(G_n\right)$ does not converge to $e_0$ on $A$. By Lemma \ref{L: general convergence almost implies unique limit}, $\left(G_n\right)$ diverges on $A$.
\end{proof}

Proposition~\ref{P: main connection} is significant because it shows that sets of divergence can be studied without M\"obius transformations. In \S\ref{S: conical} we investigate the structure of sets of divergence without using M\"obius transformations.


\section{Applications to continued fractions}\label{S: continued fractions}{\quad}\vspace{6pt}


Attached to the concept of continued fractions are particular sequences of M\"obius maps, which we describe as \emph{continued fraction sequences}. We will study both the limit sets and the conical limit sets of continued fraction sequences.

An \emph{infinite complex continued fraction} $\mathbf{K}(a_n|\, b_n)$ is a formal expression
\begin{equation}\label{E: cf}
                     \cfrac{a_1}{b_1+
                      \cfrac{a_2}{b_2+
                       \cfrac{a_3}{b_3+ \dotsb
}}}\;,
\end{equation}
where the $a_i$ and $b_j$ are complex numbers and no $a_i$ is equal to $0$. We define M\"obius transformations \begin{equation}\label{E: tns}
t_n(z)=a_n/(b_n+z),\quad \text{for } n=1,2,\dotsc\, .
\end{equation}
 We write $T_n=t_1 \circ \dotsb \circ
t_n$ for the inner composition of $t_1,\dots,t_n$. We say that the continued fraction \emph{converges classically} if the sequence $\left(T_n(0)\right)$ converges. Each M\"obius transformation $t_n$ maps $\infty$ to $0$, therefore $T_n(\infty)=T_{n-1}(0)$ for $n=2,3,\dotsc$. We define a \emph{continued fraction sequence} in $m$ dimensions to be a sequence $\left(T_n\right)$ of orientation-preserving maps in $\text{Isom}(\mathbb{H}^{m+1})$ for which  $T_n(\infty)=T_{n-1}(0)$ for every $n$, with the convention that $T_0$ is the identity transformation. This idea of generalising continued fractions to $m$-dimensional was first proposed by Beardon in \cite{Be1}. The two-dimensional concept of classical convergence extends to all dimensions. Proposition~\ref{P: two notions of general convergence} applied with $u_n=\infty$  and $v_n=0$ shows that classically convergent continued fractions converge generally. 

\subsection{Conical limit sets of continued fractions}\label{SS: CF conical sets}{\quad}\vspace{6pt}

In Proposition~\ref{P: main connection} it was shown that each conical limit set arises as the set of divergence of a sequence of M\"obius transformations. In Proposition~\ref{P: CFSDiffConical} we will show that each conical limit set arises as the set of divergence of a continued fraction sequence. Then in Proposition~\ref{P: CFconv} we will show that the class of conical limit sets that do not contain the points $0$ or $\infty$ is equal to the class of conical limit sets of classically convergent continued fractions. 

To prove these results we move to the ball model of hyperbolic space. In this model, by a continued fraction sequence we mean a sequence $\left(T_n\right)$ of hyperbolic isometries of $\mathbb{B}^{m+1}$ satisfying $T_n(e)=T_{n-1}(-e)$, where $T_0$ is the identity map and $e = (1,0, \dots,0)$. Classical convergence of $\left(T_n\right)$ means convergence of $\left(T_n(e)\right)$ in the spherical metric.  

A preliminary lemma is used in the proofs of Propositions~\ref{P: CFSDiffConical} and~\ref{P: CFconv}. In the lemma we use the notation $\gamma$ for the oriented hyperbolic geodesic from $-e$ to $e$ (tracing out a Euclidean diameter of $\mathbb{B}^{m+1}$). 

\begin{lemma}\label{L: orthogonal maps}
If $z\in\mathbb{B}^{m+1}$, $x,y\in \mathbb{S}^{m}$ and $|x-y|=2/\cosh\rho(z,\gamma)$, then there is an orientation-preserving M\"{o}bius map $U$ preserving $\mathbb{B}^{m+1}$ such that $U(z)=\jj$, $U(e)=x$ and $U(-e)=y$, provided that $m \ge 2$. For $m=1$ the map exists if $z$ is on the left (respectively right) of $\gamma$ and $0$ is on the left (respectively right) of the oriented geodesic from $y$ to $x$.
\end{lemma}
\begin{proof}
Since $\textup{Isom}(H^{m+1})$ acts transitively on $H^{m+1}$, we can choose $V \in \textup{Isom}\left(\mathbb{B}^{m+1}\right)$ with $V^{-1}(\jj)=z$. From equation~(\ref{E: key geometry}) we have
\[
|x-y| \,=\, 2/\cosh \rho(z, \gamma) \,=\, 2/\cosh \rho (\jj,V(\gamma)) \,=\, |V(e)-V(-e)|\,,
\]
So we can choose an orthogonal linear map $W$ with $W(V(e))=x$ and $W(V(-e))=y$. We can take $W$ to be orientation-preserving if $m \ge 2$, but for $m=1$ this extra condition can be satisfied if and only if $z$ is on the same side of $\gamma$ as $0$ is of $[y,x]$. Finally, define $U=W \circ V$ to obtain the required map.
\end{proof}

\begin{proposition}\label{P: CFSDiffConical}
For any conical limit set $A$ in $\mathbb{S}^m$ there is a continued fraction sequence $\left(T_n\right)$ with $A=\Lambda_c(T_n^{-1}(\jj))$.
\end{proposition}
\begin{proof} 
 Choose a sequence $(z_n)$ in $\mathbb{B}^{m+1}$ with $C=\Lambda_c(z_n)$. Using Lemma~\ref{L: orthogonal maps} we may recursively define a sequence $\left(T_n\right)$ of M\"obius transformations satisfying $T_n(z_n)=\jj$ and $T_n(e)=T_{n-1}(-e)$. 
\end{proof}

The conical limit set of a classically convergent continued fraction sequence cannot contain either of the points $0$ or $\infty$ by Proposition~\ref{P: Aebischer}. However, this is the only constraint that need be placed on a conical limit set for it to be the conical limit set of a classically convergent continued fraction.

\begin{proposition}\label{P: CFconv}
Let $A$ be a conical limit set in $\mathbb{S}^m$. Then $A$ arises as the set of divergence of a classically convergent continued fraction sequence if and only if $e,-e \notin A$.
\end{proposition}
\begin{proof}
One implication has been discussed. For the converse, suppose that $C \subseteq \mathbb{S}^m$ belongs to $\textup{CL}(m)$ and does not contain $e$ or $-e$. We choose an escaping sequence $(z_n)$ in $\mathbb{B}^{m+1}$ with $\Lambda_c\left(z_n\right) = C$. Because $-e$ and $e$ do not belong to $\Lambda_c\left(z_n\right)$, the sequence $\rho(z_n, \gamma)$ contains no bounded subsequence, and therefore we can reorder $\left(z_n\right)$ so that $\rho(z_n, \gamma)$ is an increasing sequence. In the case $m=1$ we select the $z_n$ to be on the right of $\gamma$ for odd $n$ and on the left for even $n$. (If necessary we may do this by inserting extra points into the sequence lying on the horocycles through $\jj$ which are tangent to $\mathbb{S}^m$ at $e$ and $-e$; doing this will not enlarge $\Lambda_c\left(z_n\right)$).

Define $T_0$ to be the identity map in $\text{Isom}\left(\mathbb{B}^{m+1}\right)$ and define $\theta_0=0$. Then use Lemma~\ref{L: orthogonal maps} to define inductively for each $n \ge 1$ an orientation-preserving map $T_n$ in $\text{Isom}\left(\mathbb{B}^{m+1}\right)$ such that 
\[
T_{n}(e) = T_{n-1}(-e) = \left(\cos \theta_n, \sin \theta_n, 0,\dots,0\right)\,,
\]
 where 
\[ 
\frac{\theta_n - \theta_{n-1}}{2}\, = \,  (-1)^n\, \sin^{-1}\left(\frac{1}{\cosh \rho\left(z_n, \gamma\right)}\right)\,,
\] 
and 
\[
T_n\left(z_n\right) = \jj\,.
\] 
Here we use the value of $\sin^{-1}$ in $[0, \pi/2]$. The sequence $\left(\theta_n\right)$ converges since the differences $\theta_n - \theta_{n-1}$ are alternating in sign and decreasing in modulus. Therefore $\left(T_n(e)\right)$ converges. 
\end{proof}

\subsection{Limit sets of continued fractions}\label{SS: CF limit sets}{\quad}\vspace{6pt}

The \emph{limit set} of a continued fraction is the set of accumulation points in $S^{m-1}$ of $\left(T_n^{-1}(\jj)\right)$.  In \cite{Aeb} and \cite{Be1}  it is proven that the complement of the limit set of a M\"obius sequence $\left(G_n\right)$ that is generally convergent to $x$ is the largest open set on which $\left(G_n\right)$ converges locally uniformly to $x$. That observation is inspired by similar results in discrete group theory and iteration theory. In contrast to the complicated restrictions on the structure of conical limit sets of continued fractions, the limit sets of continued fractions are not restricted in any way other than being closed. This is so even if we confine our attention to classically convergent continued fractions.

\begin{proposition}\label{P: limit sets of classically convergent CFs}
Each closed set $C$ in $\mathbb{S}^{m}$ is the limit set of a classically convergent 
continued fraction. 
\end{proposition} 
\begin{proof}
We use a similar proof to that of Proposition~\ref{P: CFconv}.The only change required is that we select an escaping sequence $\left(z_n\right)$ of points in $\mathbb{B}^{m+1}$ whose limit set is $C$ and for which $\rho(z_n, \gamma)$ is an unbounded increasing sequence, so that $e$ and $-e$ do not belong to $\Lambda_c\left(z_n\right)$. This is easy to do: simply select a sequence $\left(\zeta_n\right)$ of points in $C \setminus \{e, -e\}$ so that the accumulation set of $\left(\zeta_n\right)$ is $C$, then inductively define $z_n = r_n \zeta_n$ with $r_n \in [0,1)$ so that $r_n \to 1$ as $n \to \infty$ and $\rho(z_n, \gamma) \nearrow \infty$ as $n \to \infty$. 
\end{proof}

\section{Proofs of the Main Theorems}\label{S: conical}

Theorem~\ref{T: main 1} and Theorem ~\ref{T: main 3} were proved in the case $m=2$ in \cite{EP} and \cite{PT}, in terms of M\"{o}bius maps and sets of divergence. Our proofs are essentially geometric interpretations of their proofs, written in such a way that the generalisation to higher dimensions involves no extra work. Some extra work is required for the case $m=1$.

\subsection{Basic results on conical limit sets}\label{SS: basic}{\quad}\vspace{6pt}

We first establish some properties of conical limit sets, some of which  are used in the proof of Theorem~\ref{T: main 1}. The first of these is trivial, and may be found in \cite{PT} for the case $m=2$. Let $E$ and $F$ be subsets of $H^{m+1}$. Since any infinite sequence from $E \cup F$ contains an infinite subsequence from one of $E$ and $F$, we have $\Lambda_c(E \cup F) = \Lambda_c(E) \cup \Lambda_c(F)$.

\begin{lemma}\label{L: unions}
The class $\textup{CL}(m)$ is closed under finite unions.
\end{lemma}

The class $\textup{CL}(m)$ is not closed under countable unions, for any $m$. Indeed, any singleton is a conical limit set, but we will see in \S\ref{S: countable conical limit sets} that countable dense subsets of $S^m$ are not conical limit sets. Also, $\textup{CL(m)}$ is not closed under taking complements. For example, the set $\mathbb{Q}$ of rational numbers is not in $\textup{CL}(1)$ but $\mathbb{R} \setminus \mathbb{Q}$ is in $\textup{CL}(1)$ because it is a $G_\delta$ set. It is not known whether $\textup{CL}(m)$ is closed under finite intersections.

\begin{lemma}\label{L: intersectOpens}
The class $\textup{CL}(m)$ is closed under intersection with open sets.
\end{lemma}
\begin{proof}
By M\"{o}bius-invariance, it suffices to work in the upper half-space model and prove that if $X \subseteq \mathbb{R}^m$ belongs to $\textup{CL}(m)$ and $U \subseteq \mathbb{R}^m$ is open, then $X \cap U \in \textup{CL}(m)$.
Define  the \emph{dune} of $U$ to be the following subset of $\mathbb{H}^{m+1}$:
\[ \mathcal{D}(U) = \left\{ (t,v) \in \mathbb{H}^{m+1}\, |\,t < 1, \, t < d\left(v, \mathbb{R}^m \setminus U\right)^2\, \right\}\,.\]
Here, $d$ refers to Euclidean distance in $\mathbb{R}^m$. Choose a subset $E$ of $\mathbb{H}^{m+1}$ such that $X = \Lambda_c(E)$. Then we claim that \[\Lambda_c(E \cap \mathcal{D}(U)) = X \cap U\,.\]
Indeed, suppose $\left(w_n\right)$ is a sequence of points of $E$ that converges conically to a point $z \in \mathbb{R}^m$. If $z \in U$ then $w_n \in \mathcal{D}(U)$ for $n$ sufficiently large, but if $z \not\in U$ then $w_n \not\in \mathcal{D}(U)$ for $n$ sufficiently large.
\end{proof}

\begin{corollary}\label{C: opens}
The class $\textup{CL}(m)$ contains all open sets.
\end{corollary}
\begin{proof}
Since $S^m = \Lambda_c(H^{m+1})$, we have $S^m \in \textup{CL}(m)$, and the result follows from Lemma~\ref{L: intersectOpens}.
\end{proof}

The class $\textup{CL}(m)$ also contains all closed sets, because if $X$ is a closed subset of $S^m$ and $E=\bigcup_{x\in X}[\jj,x]$, then $X=\Lambda_c(E)$.

\subsection{Proof of Theorem~\ref{T: main 2}, part (i)}\label{SS: local property}{\quad}\vspace{6pt}

In this section we prove that being a conical limit set is a local property, in the strong sense that in order to check that $X \in \textup{CL}(m)$, one only needs to check that $X$ agrees with a conical limit set on a neighbourhood of each point of $X$. As a corollary we find that the same statement applies to the complements of conical limit sets.
 
Let $X \subseteq S^m$ and let $\mathcal{U} = \{U_i: i \in I\}$ be a cover of $X$ by open sets in $S^m$, for some indexing set $I$. Lemma~\ref{L: intersectOpens} shows that if $X \in \textup{CL}(m)$ then  $X \cap U_i  \in \textup{CL}(m)$ for each $i \in I$.  We have to prove the converse: if $U_i\cap X$ is in $\text{CL}(m)$ for all $i \in I$, then $X\in\text{CL}(m)$. 
\begin{proof}
 If $\mathcal{U}$ has a finite subcover, then the condition is sufficient, by Lemma~\ref{L: unions}. However, $X$ does not have to be compact. Because $\bigcup_{i\in I} U_i$ is an open set in $S^m$, it is a paracompact metric space. So there is a locally finite open cover $\mathcal{V}$ of $X$ subordinate to $\mathcal{U}$. Lemma~\ref{L: intersectOpens} ensures that if $V\in\mathcal{V}$ and $U\in\mathcal{U}$ satisfy $V\subseteq U$, then $X\cap V = X\cap U \cap V$ lies in $\text{CL}(m)$. Now $\mathcal{V}$ satisfies the same hypotheses as $\mathcal{U}$, and it is a locally finite open cover of $X$. So we may assume that  $\mathcal{U}$ is locally finite (and therefore countable).

We use the upper half-space model: it suffices to deal with the case where $X$ and all $U_i$ are  subsets of $\mathbb{R}^m$. For each $i\in I$, choose $E_i \subseteq \mathbb{H}^{m+1}$ with $X \cap U_i = \Lambda_c(E_i)$. Recall the definition of the dune $\mathcal{D}$ from the proof of Lemma~\ref{L: intersectOpens}.
Define 
\[
F_i = E_i \cap \mathcal{D}(U_i)\,.
\]
We will show that $X=\Lambda_c\left(\bigcup_{i \in I} F_i\right)$. Any $x \in X$ belongs to some  $U_i$, so $x \in \Lambda_c(E_{i})$ and therefore $x \in \Lambda_c(F_i)$.
 On the other hand, if $x \in \Lambda_c\left(\bigcup_{i \in I} F_i\right)$, then for some $\alpha > 0$, $x$ is the $\alpha$-conical limit of some sequence $\left((t_n, v_n)\right)$ of points in $\bigcup_{i \in I} F_i$ (where $v_n\in\mathbb{R}^m$ and $t_n>0$). We claim that there is an infinite subsequence of points all belonging to one $F_i$, and therefore $x \in \Lambda_c(F_i) \subseteq X$. Indeed, there are only finitely many sets in $\mathcal{U}$ that contain $x$, so if there is no such infinite subsequence then for each large enough $n$ there exists $i \in I$ such that $(t_n,v_n) \in F_i$ but $x \not\in U_{i}$, and therefore $t_n < d(v_n,x)^2$. Therefore $\left((t_n, v_n)\right)$ does not converge conically to $x$, which is the desired contradiction. 
\end{proof}
\begin{corollary} If $\{U_i: i \in I\}$ is an open cover of $S^m \setminus X$ such that $U_i \cap X \in \textup{CL}(m)$ for all $i \in I$, then $X \in \textup{CL}(m)$.
\end{corollary}
\begin{proof}
The set $S^m \setminus \bigcup_{i\in I} U_i$ is a member of $\textup{CL}(m)$ because it is closed. The set $\bigcup_{i\in I} (U_i \cap X)$ is a member of $\textup{CL}(m)$ by Theorem~\ref{T: main 2}(i). Also, $X$ is a member of $\textup{CL}(m)$ because $X= \left(S^m \setminus \bigcup_{i\in I} U_i\right) \cup \left(\bigcup_{i\in I} (U_i \cap X)\right)$.
\end{proof}

\subsection{Proof of Theorem~\ref{T: main 1}, parts (i) and (ii)}{\quad}

In this section we prove that the class $\textup{CL}(m)$ lies strictly between $G_\delta$ and $G_{\delta\sigma}$. 

 First we show that $\textup{CL}(m)$ consists of $G_{\delta\sigma}$ sets. This is immediate from Corollary~\ref{C: conical limit set formulae}, since each set $\bigcup_{n=p}^\infty \textup{Shad}_{\jj}^\alpha\left(w_n\right)$ is a union of open balls, therefore it is open, and so $A$ is a countable union of countable intersections of open sets.

For each $m \ge 1$ there are $G_{\delta\sigma}$ sets in $S^m$ that do not belong to $\textup{CL}(m)$. For example, we will see from Lemma~\ref{L: E is uncountable} that countable dense subsets of $S^m$ are not conical limit sets, yet they are $G_{\delta\sigma}$ sets.

We move on to Theorem~\ref{T: main 1} (ii). Let $B(c,r)$ denote the open ball in $\mathbb{R}^m$ of radius $r>0$, centred on $c\in\mathbb{R}^m$.
\begin{lemma}\label{L: covers of finite multiplicity}
Suppose that $V$ is an open subset of $\mathbb{R}^m$, and let $\epsilon >0$. Then there exists a sequence of open balls $B_i = B_i(v_i,t_i)$ such that
\begin{enumerate}
\item $V = \bigcup_{i=1}^\infty B_i$; 
\item $t_i \leq \epsilon$;
\item for each $\alpha>0$, each point in $\mathbb{R}^m$ lies in only finitely many of the balls $B(v_i,\alpha t_i)$; \item $t_i/d(v_i,\partial V) \to 0$ as $i \to \infty$.\end{enumerate} 
\end{lemma}
\begin{proof}
 Define 
\[
W_0 = \{x \in \mathbb{R}^m\,:\, d(x, \mathbb{R}^m \setminus V) \geq \epsilon\}
\] and, for $n=1,2, \dots$, 
\[
W_n = \{x \in \mathbb{R}^m \, : \, \epsilon 2^{-n} \leq d(x, \mathbb{R}^m \setminus V) \leq \epsilon 2^{1-n}\}\,.
\]
 Each $W_n$ is a closed subset of $V$, therefore each $W_n$ is compact. So $W_n$ may be covered by finitely many open balls of the form $B\left(v, \epsilon 2^{-2n}\right)$, $(v \in W_n)$. Call these balls $B_{nj}$, $j=1, 2,  \dots, k(n)$. We claim that the countable collection of all such balls has the required properties (we can order these balls in a sequence in an arbitrary fashion). First notice that $V$ is the union of the $W_n$ so $V$ is covered by the $B_{nj}$ (condition (i)). Condition (ii) holds as $\epsilon 2^{-2n}\leq \epsilon$. Next, given $\alpha>0$ we may check that for sufficiently large $n$ the balls $B_{nj}(v_{nj},\alpha\epsilon 2^{-2n})$ are contained within $W_{n-1}\cup W_n \cup W_{n+1}$, hence condition (iii) is satisfied. Last, $t_{nj}/d(v_{nj},\partial V) \leq \epsilon 2^{-2n}/\epsilon 2^{-n}=2^{-n}$, so condition (iv) is also satisfied.
\end{proof}

We can now prove that every $G_\delta$ subset of $S^m$ is a conical limit set. This generalises \cite[Theorem 3]{PT} to arbitrary dimensions. In the proof, we make use of the observation that M\"obius transformations preserve the class $\text{CL}(m)$.

\begin{proof}
We work in the upper half-space model of hyperbolic space. It suffices to show that any $G_\delta$ subset $E$ of $\mathbb{R}^m$ is a conical limit set, since $\mathbb{R}^m_\infty\in\text{CL}(m)$ and all other $G_\delta$ sets can be mapped into $\mathbb{R}^m$ by a M\"obius map.
Let
 \[
E = \bigcap_{n=1}^{\infty} U_n\,,
\] 
where each set $U_n\subseteq\mathbb{R}^m$ is an open subset of $\mathbb{R}^m$. Define $V_n = \bigcap_{m=1}^n U_n$ so that $V_n$ is a decreasing sequence of open sets whose intersection is $E$. Cover $V_n$ by a sequence of balls $B_{nj}$, $j=1,2, \dots$ satisfying the conditions of Lemma~\ref{L: covers of finite multiplicity} with $\epsilon = 1/n$. Rearrange all these balls into a simple sequence $B_1, B_2, \dots$, where $B_n = B\left(v_n, t_n\right)$. Let $w_n$ be the point $(t_n,v_n)$ in $\mathbb{H}^{m+1}$ for $n=1,2,\dotsc$, so that  $B_n=\textup{Shad}_\infty^\alpha(w_n)$, where $\sinh\alpha = 1$ by Lemma~\ref{L: alpha shadows from infinity}. Lemma~\ref{L: shadow descriptions of conical limit sets} (ii) now shows that $E \subseteq \Lambda_c^\alpha\left(w_n\right)$ for any $\alpha > \sinh^{-1}(1)$. 

We must also check that $\Lambda_c(w_n) \subseteq E$. We show that if $z \not\in E$ then for each $\alpha > 0$, the point $z$ lies in only finitely many of the balls $B\left(v_i, \alpha t_i\right)$. Choose $N$ large enough that for $n \ge N$, we have $z\notin V_n$ and $t_n/d(v_n,\partial V_n) < 1/\alpha$. For each $n< N$, the point $z$ lies in only finitely many $B_{nj}(v_{nj},\alpha t_{nj})$, by Lemma~\ref{L: covers of finite multiplicity}. For each $n>N$ we have that
\[ 
d(v_{nj},z)\geq d(v_{nj},\partial V_n) > t_{nj}\alpha,
\]
and therefore $z\notin B(v_{nj},\alpha t_{nj})$. So for each $\alpha > 0$, $z$ lies in $\textup{Shad}_{\infty}^\alpha(w_n)$ for only finitely many $n$.
This completes the proof.
\end{proof}

To finish proving Theorem~\ref{T: main 1}(ii) we must show that the class of $G_\delta$ sets in $S^m$ is not equal to the class $\text{CL}(m)$.  For $m \ge 2$ this will follow from Theorem~\ref{T: main 1}(iii), which says that all $G_{\delta\sigma}$ subsets of a codimension-one sphere in $S^m$ belong to $\textup{CL}(m)$; for $m \ge 2$ not every such set is a $G_\delta$ set. Notice that for a subset $A \subseteq S^{m-1} \subseteq S^m$ we may have $A \in \textup{CL}(m)$ but $A \not\in \textup{CL}(m-1)$. In contrast, the Borel hierarchy is unchanged when we consider subsets of $S^{m-1}$ as being subsets of $S^m$. For $m=1$ we give instead an explicit construction of a one-dimensional conical limit set that is not a $G_\delta$ set.

\begin{figure}[ht]
\centering
\includegraphics[scale=0.7]{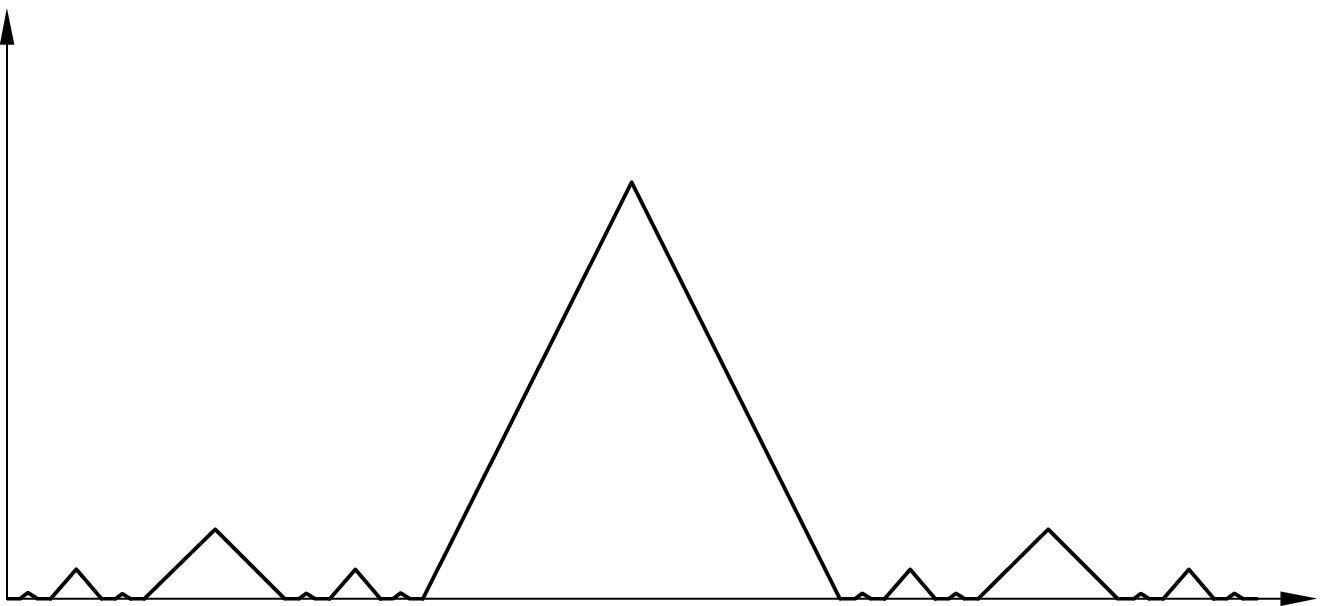}
\caption{}
\label{F: graph}
\end{figure}

 Consider the middle-thirds Cantor set $C=\bigcap_{n=1}^\infty J_n$, where $J_1=[0,1]$, $J_2 =[0,1/3]\cup [2/3,1]$, $J_3 =[0,1/9]\cup [2/9,1/3]\cup [2/3,7/9]\cup [8/9,1]$, and so on. Let $A$ be the countable set $\bigcup_{n=1}^\infty \partial J_n$ (the accessible points of $C$). Since $A$ is contained in its own derived set, if it were a $G_\delta$ set then it would be uncountable. (This follows from a standard argument, and we will use a similar argument in the proof of Lemma~\ref{L: E is uncountable}.) As $A$ is countable, it is not a $G_\delta$ set. However, $A \in \textup{CL}(1)$. To see this, let $E$ be the graph (shown in Figure~\ref{F: graph}) of the function $f: [0,1] \setminus C \to \mathbb{R}^+$ defined by $f(x) = d(x,C)/n$ when $x$ lies in an open interval of $[0,1] \setminus C$ of length $3^{-n}$. Then it is easy to check that $A = \Lambda_c(E)$.  Theorem~\ref{T: main 3} gives us another way to check whether a given countable set is a conical limit set.
\subsection{Proof of Theorem~\ref{T: main 1}, part (iii)}\label{SS: hyperplanes}
\begin{figure}[hb]
\centering
\includegraphics[scale=0.5]{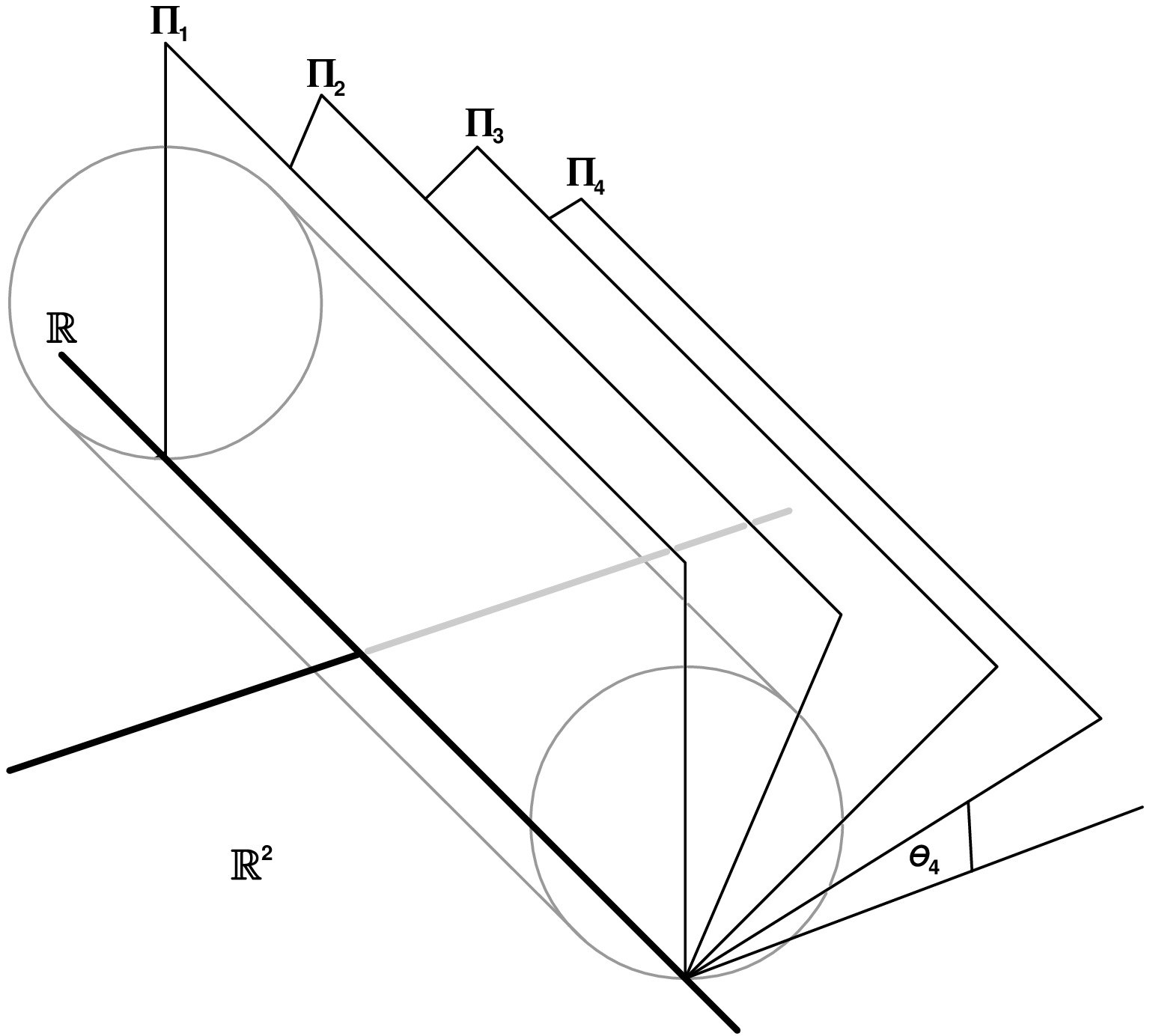}
\caption{}
\label{F: planes}
\end{figure}

{\quad}\vspace{6pt}\\
Consider a $G_{\delta\sigma}$ set $A$ in $\mathbb{R}$. We have seen that $A$ is not necessarily a member of $\textup{CL}(1)$. However, the real line may be embedded as the real axis of $\mathbb{R}^2$ and we shall see that it \emph{is} true that $A\in\text{CL}(2)$ when we consider $A$ as a subset of $\mathbb{R}^2$.   More generally, let $Y$ be a sphere of codimension one in $S^m$. We will show that  all $G_{\delta\sigma}$ subsets of $Y$ belong to $\textup{CL}(m)$.
\begin{proof}
Inside $\mathbb{H}^{m+1} = \{ (x_0, x_1, \dots, x_m) \, | \, x_0 > 0 \}$, we have $\mathbb{H}^m$ isometrically embedded as the hyperplane $x_m = 0$, and this extends to an embedding of the ideal boundaries $\mathbb{R}^{m-1}_{\infty} \subseteq\mathbb{R}^m_\infty$. By M\"{o}bius-invariance of conical limit sets, it suffices to prove that any $G_{\delta\sigma}$ subset $X \subseteq \mathbb{R}^{m-1}_\infty$ is the conical limit set of some subset of $\mathbb{H}^{m+1}$, and as usual it suffices to consider $X \subseteq\mathbb{R}^{m-1}$. Express $X$ as $\bigcup_{n=1}^{\infty} X_n$ for $G_\delta$ sets $X_n \subseteq\mathbb{R}^{m-1}$. By Theorem~\ref{T: main 1}(ii), $X_n = \Lambda_c(\Gamma_n)$ for some subset $\Gamma_n \subseteq \mathbb{H}^m$. Note that $\Lambda_c\left(\bigcup_{n=1}^\infty \Gamma_n\right)$  may contain points that do not belong to any $X_n$. For example, some  sequence $\left(w_n\right)$, where  $w_n \in \Gamma_n$, may converge conically to a point not in $X$. We therefore define sets 
\[ \Gamma_n' = \left\{ \left(x_0\sin\tfrac{\pi}{n}, x_1, x_2, \dots, x_{m-1}, x_0\cos\tfrac{\pi}{n}\right)\,:\, \left(x_0, \dots, x_{m-1}, 0\right) \in  \Gamma_n, \,x_0 < 2 \sin\tfrac{\pi}{n}\,\right\}\,. \]
 We obtain $\Gamma_n'$ by first rotating $\Gamma_n$ about $\mathbb{R}^{m-1}$ and then truncating the rotated sets $\Gamma_n$ so that all points are contained in a Euclidean cylinder of radius $1$, tangent to $\mathbb{R}^{m}$ along $\mathbb{R}^{m-1}$.   See Figure~\ref{F: planes}. Any sequence from $\bigcup_{n=1}^\infty \Gamma_n'$ that converges conically to a point of $\mathbb{R}^{m}$ must therefore converge to a point of $\mathbb{R}^{m-1}$, and it may only visit finitely many of the planes $\Pi_n$. Therefore its limit is already contained in some $X_n$. Conversely, for each conically convergent standard sequence $(z_n)$ in $\Gamma_n$ (standard sequences were defined in Definition~\ref{D: standard and escaping}), there is a corresponding infinite subsequence in $\Gamma_n'$ which converges conically to the same limit point. Therefore 
 \[\Lambda_c\left(\bigcup_{n=1}^\infty \Gamma_n'\right) \,=\, \bigcup_{n=1}^\infty \Lambda_c(\Gamma_n') \, = \, \bigcup_{n=1}^\infty X_n \, = \, X\,.\] 
\end{proof}
\subsection{Proof of Theorem~\ref{T: main 3}}\label{S: countable conical limit sets}{\quad}\vspace{6pt}

Erd\H{o}s and Piranian \cite[Theorem 2]{EP} gave a characterisation of countable sets of divergence in $\Chat$. Our Theorem~\ref{T: main 3} reinterprets their theorem in terms of conical limit sets and hyperbolic geometry, and generalises it to arbitrary dimensions. The hyperbolic viewpoint does not simplify the proof much, but we believe it gives an interesting way of understanding the theorem. 

For any subset $E \subseteq\mathbb{C}$, Erd\H{o}s and Piranian defined a subset of \emph{good points} of $\mathbb{C}$ with respect to $E$. Here is their definition, rewritten to allow arbitrary dimension.
\begin{definition}\label{D: gd(E)}
Let $m \geq 1$ and  $E \subseteq \mathbb{R}^m$ be given. Then $z \in \textup{gd}(E)$ if and only if the following holds. For each $\epsilon>0$ there exists $\delta>0$ such that, whenever $0<|v-z| < \delta$, there exists $x \in E$ such that $|x-v|\, < \,\epsilon\, |v-z|$.
\end{definition} 
The next lemma describes $\mathbb{R}^m \setminus \textup{gd}(E)$ in terms of shadow sets, using Lemma~\ref{L: alpha shadows from infinity}. We use the notation $\delta_z$ for the `vertical' geodesic in $\mathbb{H}^{m+1}$ between $z$ and $\infty$. We use the notation $B(z,r)$ for the open Euclidean ball of radius $r$, centred on $z$.
\begin{lemma}\label{L: hyperbolic definition of gd}
A point $z$ lies in $\mathbb{R}^m \setminus \textup{gd}(E)$ if and only if there exist $\alpha > 0$ and a sequence $(w_n)$ in $\mathbb{H}^{m+1}$ such that $w_n \to z$ ideally and for each $n$, $\rho\left(w_n, \delta_z\right) \leq 2\alpha$ and $\textup{Shad}_\infty^{\alpha}\left(w_n\right) \cap E = \emptyset$. Moreover, in this case such a sequence $\left(w_n\right)$ exists for arbitrarily large $\alpha$.
\end{lemma}
\begin{proof}
First observe that $z\in \mathbb{R}^m \setminus \text{gd}(E)$ if and only if there exists an $\epsilon>0$ and a sequence $\left(v_n\right)$ in $\mathbb{R}^m$ such that $v_n\rightarrow z$ and the balls $B(v_n,\epsilon |v_n-z|)$ are all disjoint from $E$.

 Suppose that $\left(w_n\right)$ and $\alpha$ exist as in the lemma. Writing $w_n = (t_n, v_n) \in \mathbb{H}^{m+1}$, we have $|v_n - z| \to 0$ as $n \to \infty$, but $|v_n - z| \leq \sinh (2 \alpha) t_n$ for all $n$. Let $\epsilon = 1/(2\cosh\alpha)$. If $x\in B(v_n,\epsilon|v_n-z|)$  then $|x - v_n| < (\sinh \alpha) t_n$, that is, $x \in \textup{Shad}_\infty^\alpha\left(w_n\right)$, so $x \not\in E$. Hence $z \not \in \textup{gd}(E)$.

Conversely, suppose $z \not\in \textup{gd}(E)$. Then there exists $\epsilon > 0$ (which we can choose to be arbitrarily small) and a sequence of points $\left(v_n\right)$ in $\mathbb{R}^m$ such that $v_n \to z$ as $n \to \infty$ and \[( |x-v_n| < \epsilon |v_n - z| )\, \implies\, x \not\in E\,.\] Take $\alpha = \cosh^{-1}(1/(2\epsilon))$ and set $t_n =  |v_n - z|/\sinh(2\alpha)$, and $w_n = (t_n, v_n)$. Then $\rho\left(w_n, \delta_z\right) = 2\alpha$. If $x \in \textup{Shad}_\infty^{\alpha}\left(w_n\right)$ then \[|x-v_n| < (\sinh \alpha)t_n = |v_n - z|/(2 \cosh \alpha) = \epsilon |v_n - z|\,,\] so $x \not\in E$. Thus $\textup{Shad}_\infty^{\alpha}\left(w_n\right) \cap E = \emptyset$.
\end{proof}

We can now give a characterization of $\textup{gd}(E)$ using hyperbolic geometry.
\begin{lemma}
For any subset $E$ of $\mathbb{R}^m$ we have 
\[
 \textup{gd}(E) =  \mathbb{R}^m \setminus  \Lambda_c( \mathbb{H}^{m+1} \setminus \overline{\textup{co}}(E)),
\]
 where $\overline{\textup{co}}(E)$ is the hyperbolic convex hull of $E$ in $\mathbb{H}^{m+1}$. Thus 
\[
E \cap \textup{gd}(E) = E \setminus \Lambda_c(\mathbb{H}^{m+1} \setminus \overline{\textup{co}}(E)).
\]
\end{lemma}
\begin{proof}
 A point $w \in \mathbb{H}^{m+1}$ is not in $\overline{\textup{co}}(E)$ if and only if there exists a hyperbolic hyperplane $P$ separating $w$ from $E$. If $(w_n)$ converges $\alpha$-conically to $x \in \mathbb{R}^m$ and each $w_n \not\in \overline{\textup{co}}(E)$ then we can choose corresponding hyperbolic hyperplanes $P_n$, so that $P_n$ separates $w_n$ from $E$. The $(m-1)$-dimensional sphere or hyperplane in which $P_n$ meets $\mathbb{R}^{m}$ separates $\mathbb{R}^m$ into two components. One contains $E$ and the other contains a ball $B(a_n,r_n) \subseteq\mathbb{R}^m$ such that $r_n\geq k|a_n-x|$ and $|a_n-x| < 2|w_n - x|$, for some $k>0$. This can easily be verified by geometric means. It follows that $x \not\in \textup{gd}(E)$. For the converse, suppose that $x \not\in \textup{gd}(E)$. Then we find a sequence of open balls in $\mathbb{R}^m$ that converges to $x$, and each ball subtends an angle at $x$ greater than or equal to some $\theta > 0$. Each such ball is the ideal boundary of a hyperbolic half-space in $\mathbb{H}^{m+1}$ whose distance $\alpha$ from the geodesic joining $x$ and $\infty$ depends only on $\theta$. We can then find a sequence of points, one in each of these half-spaces, that converges $(1+\alpha)$-conically to $x$. 
\end{proof}
 
Given $E\subseteq \mathbb{R}^m$, we define inductively a set $E^\chi$ for each ordinal $\chi$, as follows:
\begin{eqnarray*} 
 E^1 & = & E,  \\
 E^{\chi + 1} & =&  E^\chi \cap \textup{gd}(E^\chi) 
 \, , \quad\text{and}\\
 E^{\chi} & = & \bigcap_{\psi < \chi} E^\psi, \quad \text{for a limit ordinal $\chi$}.
\end{eqnarray*}  

Theorem~\ref{T: main 3} states that a countable set $E \subseteq S^m$ is in $\textup{CL}(m)$ if and only if there exists an ordinal $\chi$ such that $E^{\chi} = \emptyset$. Let $E$ be a countable subset of $\mathbb{R}^m_\infty$. Since $\psi < \chi$ implies that $E^{\chi} \subseteq E^{\psi}$, there must exist a countable ordinal $\chi$ for which $E^{\chi}=E^{\chi +1}$: otherwise there would be an uncountable descending chain of subsets of $E$, and $E$ would therefore be uncountable. In Lemma~\ref{L: E is uncountable} we show that if $E$ is also a conical limit set then $E^\chi$ must be empty, thereby establishing one implication in Theorem~\ref{T: main 3}. The converse implication is then established at the end of this section.

Before we begin the proof of Theorem~\ref{T: main 3}, we point out two simple corollaries of that theorem.
\begin{corollary}
Every countable member of $\textup{CL}(m)$ is nowhere dense. 
\end{corollary}
\begin{proof}
If a set $E$ is dense in a ball $B$ then $\text{gd}(E)\supseteq B$, hence $E^\chi \supseteq E \cap B$ for all $\chi$. 
\end{proof}
\begin{corollary}
If $E \in \textup{CL}(m)$ is countable, then every subset of $E$ is also in $\textup{CL}(m)$.
\end{corollary}
\begin{proof}
Transfinite induction shows that $F \subseteq E$ implies $F^\chi \subseteq E^\chi$ for each ordinal $\chi$.
\end{proof}

Recall the notion of a standard sequence from Definition~\ref{D: standard and escaping}.

\begin{lemma}\label{L: E is uncountable}
Suppose that $E \in \textup{CL}(m)$ and that $E$ contains a non-empty subset $F$ such that $F \subseteq \textup{gd}(F)$. Then $E$ is uncountable. 
\end{lemma}
\begin{proof}
If $E=S^m$ then $E$ is indeed uncountable, so we may use the upper
half-space model and suppose that $E \subseteq \mathbb{R}^m$ and that $E =
\Lambda_c\left(w_n\right)$, where $(w_n)$ is a standard sequence. The set
$F$ is non-empty, so we may choose a point $x(1) \in F$ and a closed ball
$D(1) \subseteq\mathbb{R}^m$ that contains $x(1)$ in its interior.
Suppose that points $x(1), \dots, x(2^n -1)$ in $F$ and disjoint closed balls $D(1),
\dots , D(2^n-1)$ have already been defined for some $n \geq 1$, in such a
way that $x(i)$ is in $D(i)^\circ$, the interior of $D(i)$, for each $i$. We will now fix
$k$ such that $2^{n-1} \leq k < 2^n $ and explain how to define $x(2k)$ and
$x(2k+1)$. Since $x(k) \in E$, there exists a subsequence $(w_{n_j})$ that
converges $\alpha$-conically to $x(k)$ for some $\alpha > 0$. Pick a
shadowing sequence $(z_j)$ such that $\rho\left(z_j, w_{n_j}\right) \le 3$
but $x(k) \not\in \textup{Shad}_\infty^2\left(z_j\right)$. So \[x(k) \in
\Lambda_c^{\alpha + 3}\left(z_j\right) \setminus
\Lambda_c^{2}\left(z_j\right)\,.\] Consider the balls
$\textup{Shad}_\infty^1\left(z_j\right)$ and
$\textup{Shad}_\infty^{\alpha+3}\left(z_j\right)$ in $\mathbb{R}^m$, for each integer $j$.
They have a common centre $v_j$, and the ratio of their radii does not
depend on $j$. Now, $x(k) \in \textup{Shad}_\infty^{\alpha+3}\left(z_j\right)$
for all $j$, and $x(k) \in \textup{gd}(F)$, so for sufficiently large $j$ the
ball $\textup{Shad}_\infty^1\left(z_j\right)$ also contains a point of $F$.
We can choose two values $j$ and $j'$ large enough so that the closed balls
$D(2k) = \overline{\textup{Shad}_\infty^1\left(z_j\right)}$ and $D(2k+1) =
\overline{\textup{Shad}_\infty^1\left(z_{j'}\right)}$ are \emph{disjoint}
and both are contained in $D(k)$. Now choose $x(2k) \in F \cap
D(2k)^{\circ}$ and $x(2k+1) \in F \cap D(2k+1)^{\circ}$. Note that  \[x(2k)
\in \textup{Shad}_\infty^{1}\left(z_j\right) \subset
\textup{Shad}_\infty^{1+3}\left(w_{n_j}\right)\,,\] and similarly $x(2k+1)
\in \textup{Shad}_\infty^4\left(w_{n_{j'}}\right)$.

 In this way we define $x(i)$ and $D(i)$ inductively for all $i \in
\mathbb{N}$. Now consider the set \[K = \{ x \in \mathbb{R}^m \,:\, x \in
D(i) \,\text{for infinitely many values of $i$}\,\}\,.\] $K$ is an
uncountable Cantor set contained in $\Lambda_c^4\left(w_n\right)$, which is
a subset of $E$.
\end{proof}

In Lemma \ref{L: E is uncountable} we do not assume that the subset $F$ is closed; if we do then there is a one-line proof: $F$ is a perfect set, hence uncountable.

In the next lemma, the notation $\Lambda\left(w_n\right)$ is used for the limit set in the ideal boundary sphere of a sequence $(w_n)$ in hyperbolic space. It contains but is not in general equal to $\Lambda_c\left(w_n\right)$.
\begin{lemma}\label{L: dense ordinary implies uncountable conical}
Let $\left(w_n\right)$ be a standard sequence in $\mathbb{H}^{m+1}$. Then $\left( \Lambda\left(w_n\right) \right)^\circ \subseteq \overline{\Lambda_c\left(w_n\right)}$. 
Moreover, if $\Lambda\left(w_n\right)$ has non-empty interior then $\Lambda_c\left(w_n\right)$ is uncountable. 
\end{lemma}
\begin{proof}
The proof is very similar to that of Lemma \ref{L: E is uncountable}. Let $D(1)$ be any closed ball contained in $\Lambda\left(w_n\right)$, and pick $x_1 \in D(1)^\circ$. We will show that $D(1)$ contains uncountably many points of $\Lambda_c\left(w_n\right)$. In particular this shows that $\Lambda_c\left(w_n\right)$ is dense in $\Lambda\left(w_n\right)^\circ$.  The construction is inductive. Suppose we have already chosen a closed ball $D(k)$ and a point $x_k \in  D(k)^\circ$. We show how to choose $D(2k), D(2k+1)$ and points $x_{2k}, x_{2k+1}$ in their interiors. We choose a subsequence $\left(w_{n_j}\right)$ converging to $x_k$, not necessarily conically. Since $\Lambda(w_n) \supset D(k)$, we can choose points $x(2k), x(2k+1)$ of $\Lambda(w_n)$ in the interiors of balls $D(2k) = \overline{\textup{Shad}_\infty^1(w_{n_j})}$ and $D(2k+1) = \overline{\textup{Shad}_\infty^1(w_{n_{j'}})}$ respectively, so that $D(2k) \cup D(2k+1) \subseteq D(k)$ and $D(2k) \cap D(2k+1) = \emptyset$. Then each infinite nested sequence of the closed balls $D(i)$ contains a unique point, and it is the 1-conical limit point of a subsequence of $\left(w_n\right)$.
\end{proof}

Recall from the discussion at the beginning of \S\ref{SS: CF limit sets} that a sequence of M\"obius transformations that converges generally to a limit $x$ also converges locally uniformly to $x$ on the complement of its limit set.

\begin{corollary}
Let $\left(G_n\right)$ be a sequence of M\"{o}bius maps acting on $S^m$ that converges generally to $x$, and let $U \subseteq S^m$ be a non-empty open set. Either $G_n$ converges locally uniformly to $x$ on a dense open subset of $U$ or there are uncountably many points of $U$ at which $\left(G_n\right)$ diverges.
\end{corollary}

It remains to prove the converse implication in Theorem \ref{T: main 3}. We must show how to exhibit a given countable set $E$ as a conical limit set, when $E^\chi = \emptyset$ for some countable ordinal $\chi$.  Since $E$ is countable, we may list $E$ as $\{z_1, z_2, z_3, \dots\}$. Since $E^\chi = \emptyset$, to each point $z_k \in E$ we can associate the unique countable ordinal $\beta(k)$ such that $z_k \in E^{\beta(k)} \setminus E^{\beta(k)+1}$.  Let $\alpha_k \ge k$ be a constant that plays the role of $\alpha$ in Lemma~\ref{L: hyperbolic definition of gd}, corresponding to the statement that $z_k \not\in \textup{gd}\left(E^{\beta(k)}\right)$. 

 We will construct a doubly-indexed sequence $\left(w_{kj}\right)$ of points in $\mathbb{H}^{m+1}$, where $k$ and $j$ run over $\mathbb{N}$, such that $\Lambda_c\left(w_{kj}\right) = E$. We write $w_{kj} = \left(t_{kj}, v_{kj}\right)$ in the usual way, and define the ball $D_{kj}^{\alpha} = \textup{Shad}_\infty^\alpha(w_{kj})$. Recall from Corollary~\ref{C: conical limit set formulae} that $z \in \Lambda_c^{\alpha}\left(w_{kj}\right)$ if and only if $z$ lies in infinitely many of the balls $D_{kj}^{\alpha}$.  We further define a ball $N_{kj}$ in $\mathbb{R}^m_\infty$ for each index pair $(k,j)$, by \[N_{kj} = D_{kj}^{\alpha_k}\,.\] 
  
 Our aim is to construct $\left(w_{kj}\right)$ so as to satisfy the following four conditions:

\begin{enumerate}
\item $\left(w_{kj}\right)$ converges conically to $z_k$ as $j \to \infty$, for each fixed $k$;
\item the diameter of $N_{kj}$ tends to zero as $j+k \to \infty$;
\item if $\beta(k) = \beta(n)$ then $N_{kj} \cap N_{ni} = \emptyset$ for all $i,j$;
\item if $\beta(k) < \beta(n)$ and $N_{kj} \cap N_{ni} \neq \emptyset$, then $N_{kj} \subseteq N_{ni}$.
\end{enumerate}

Before giving the construction, we show that these four properties suffice.

  Because of conditions (iii) and (iv), the collection of balls $N_{kj}$ that contain a given point $z$ is totally ordered by inclusion. Condition (ii) shows that if this collection is non-empty then it has a maximal element. By conditions (iii) and (iv), the corresponding ordinals $\beta(k)$ form a strictly decreasing sequence, which therefore has finite length. So each point $z \in \mathbb{R}^m_{\infty}$ lies in only finitely many of the balls $N_{kj}$. 

Suppose $z \in \Lambda_c^\alpha\left(w_{kj}\right)$. Then $z$ must be a conical limit of $\{ w_{kj}\,:\, \alpha_k < \alpha\}$, since $z$ lies in only finitely many of the balls $D_{kj}^\alpha$ for which $\alpha_k \geq \alpha$. There are only finitely many values of $k$ for which $\alpha_k < \alpha$, so condition (i) implies that $z \in \{z_k \, : \, \alpha_k < \alpha\}$. In particular, $z \in E$. So $\Lambda_c\left(w_{kj}\right) \subseteq E$. On the other hand, condition (i) ensures that $E \subseteq \Lambda_c\left(w_n\right)$. Thus $\Lambda_c\left(w_{kj}\right) = E$, as required.

To construct the double sequence $w_{kj}$, fix some listing of the ordered pairs $(k,j) \in \mathbb{N} \times \mathbb{N}$, say $((k_1,j_1), (k_2, j_2), \dots)$. We construct the sequence $w_N = w_{k_N j_N}$ inductively over $N$, making sure at each step that conditions (ii)--(iv) above and also the following auxiliary conditions are satisfied:
\begin{enumerate}
\setcounter{enumi}{4}
\item  for each $k,j \in \mathbb{N}$, $N_{kj} \cap E^{\beta(k)} = \emptyset$; 
\item  $\rho\left(w_{kj}, \gamma_{z_k}\right) < 2\alpha_k + \log 2$;
\item  $E \cap \partial N_{kj} = \emptyset$. 
\end{enumerate}
Condition (vi) will ensure that condition (i) holds.

 Suppose that $w_1, \dots ,w_{N-1}$ have already been chosen, for some $N \geq 1$. Write $k=k_N$, $j=j_N$, so that we now have to choose $w_{kj}$. The definition of $\alpha_k$ using Lemma \ref{L: hyperbolic definition of gd} tells us that there exist points $w  = (t,v) \in \mathbb{H}^{m+1}$ with $|w - z_k|$ arbitrarily small, $\rho\left(w, \gamma_{z_k}\right) \leq 2 \alpha_k$ and $\textup{Shad}_\infty^{\alpha_k}(w) \cap E^{\beta(k)} = \emptyset$. For $|w - z_k|$ sufficiently small, the diameter of $\textup{Shad}_\infty^{\alpha_k}(w)$ is less than $1/kj$. Only finitely many balls $N_{ni}$ have already been defined, and $z_k$ is not on the boundary of any of them since $z_k \in E$. Thus we can take $|w-z_k|$ sufficiently small so that $\textup{Shad}_\infty^{\alpha_k}(w)$ is contained in those previously defined balls $N_{ni}$ that contain $z_k$, and disjoint from those that do not contain $z_k$. If a ball $N_{ni}$ contains $z_k$ then $N_{ni}$ meets $E^{\beta(k)}$, so $\beta(n) > \beta(k)$. Thus if we took $w_{kj}=w$ then conditions (ii)--(vi) would be satisfied, but condition (vii) could fail. Instead we define $w_{kj} = (\delta t, v)$,  where the factor $1/2 \leq \delta \leq 1$ is chosen to ensure that condition (vii) holds; a suitable value of $\delta$ exists since $E$ is countable and therefore only countably many values of $\delta$ are unsuitable. Note that $N_{kj} = \textup{Shad}_\infty^{\alpha_k}\left(w_{kj}\right) \subseteq\textup{Shad}_\infty^{\alpha_k}(w)$, and \[\rho\left(w_{kj}, \gamma_{z_k}\right) \,\leq\, \rho\left(w, \gamma_{z_k}\right) + \rho\left(w, w_{kj}\right) \, \leq\, \rho\left(w, \gamma_{z_k}\right) + \log 2 \,.\] Thus all the conditions are satisfied and we have completed the induction step.

This completes the proof of Theorem \ref{T: main 3}.


\subsection{Proof of Theorem~\ref{T: main 1}, part (iv)}\label{SS: topology}{\quad}\vspace{6pt}

In this section we prove that the class $\textup{CL}(m)$ is not topologically defined. We do this using the characterization of countable conical limit sets given in Theorem~\ref{T: main 3}. Before we begin, we note the following trivial lemma.
\begin{lemma}\label{L: gd under Cartesian product}{\quad}\\
 Suppose that $E_1 \subseteq\mathbb{R}^m$ and $E_2 \subseteq\mathbb{R}^n$. Considering $E_1 \times E_2$ as a subset of $\mathbb{R}^{m+n}$, we have $\textup{gd}\left(E_1 \times E_2\right) = \textup{gd}\left(E_1\right) \times \textup{gd}\left(E_2\right)$.
\end{lemma}

To prove part (iv) of Theorem~\ref{T: main 1}, the main task is the one-dimensional case. We will define an increasing homeomorphism $\varphi: \mathbb{R} \to \mathbb{R}$, a countable set $A \subseteq\mathbb{R}$, and its image $B = \varphi(A)$, in such a way that $A \subseteq \textup{gd}(A)$ but $B \cap \textup{gd}(B) = \emptyset$. It will then follow that $A\notin \text{CL}(1)$ and $B\in \text{CL}(1)$, as required for the case $m=1$. For the higher-dimensional cases, we observe that applying $\varphi$ co-ordinatewise yields a homeomorphism of $\mathbb{R}^m$, which sends the $m$-fold Cartesian product $A^m$ onto $B^m$. By Lemma~\ref{L: gd under Cartesian product}, $A^m \subseteq \textup{gd}\left(A^m\right)$ but $\textup{gd}\left(B^m\right) = \emptyset$, so $A^m \notin \textup{CL}(m)$ while $B^m \in \textup{CL}(m)$.  

For each $n \in \mathbb{Z} \setminus\{0\}$, we define two affine functions
\begin{eqnarray*} f_n(x) & = & \frac{1}{2n} + \frac{x}{8n^2}, \qquad \text{and} \\
                  g_n(x) & = & \frac{\textup{sgn}(n)}{4^{|n|}} + \frac{x}{8\,\cdot\,4^{|n|}}\,.
\end{eqnarray*}
These maps are chosen so that the images $f_n([-1,1])$ are disjoint closed intervals for distinct $n$, and likewise the images $g_n([-1,1])$ are disjoint for distinct $n$. Let $\mathcal{F}$ be the semigroup generated by all the maps of the form $f_n$, ($n \in \mathbb{Z}\setminus\{0\}$), and $\mathcal{G}$ be the semigroup generated by the $g_n$. Define $A = \{ f(0) \, : f \in \mathcal{F}\}$ and $B = \{ g(0) \, : \, g \in \mathcal{G}\}$.  Both $A$ and $B$ are countable sets with self-similar structures:
\[A = \{0\} \cup \bigcup_{n \in \mathbb{Z} \setminus\{0\}} f_n(A)\,,\]
\[B = \{0\} \cup \bigcup_{n \in \mathbb{Z} \setminus\{0\}} g_n(B)\,.\]
Because $1/2n \in A$ for each $n \in \mathbb{Z} \setminus \{0\}$, we have $0 \in \textup{gd}(A)$. Using the self-similarity of $A$ we deduce that $A \subseteq \textup{gd}(A)$. However, $B$ is disjoint from all intervals of the form $(2/4^{n}, 3/4^{n})$ and therefore $0 \not\in \textup{gd}(B)$. Using the self-similarity of $B$ we find that $B \cap \textup{gd}(B) = \emptyset$. It remains to construct a homeomorphism of $\mathbb{R}$ that maps $A$ onto $B$. Note that there is an obvious one-to-one order-preserving correspondence $\varphi_0$ between the points of $A$ and those of $B$, that sends $f_{n_1} f_{n_2} \dotsb f_{n_k} (0)$ to $g_{n_1} g_{n_2} \dotsb g_{n_k} (0)$. It is easy to check that this extends to an order-preserving homeomorphism $\varphi_1$ from the closure of $A$ to the closure of $B$, and that the open intervals that make up $\mathbb{R} \setminus \overline{A}$ and $\mathbb{R} \setminus \overline{B}$ are also naturally in correspondence with each other, in an order-preserving manner. We can therefore interpolate linearly on each complementary interval to extend $\varphi_1$ to a homeomorphism $\varphi: \mathbb{R} \to \mathbb{R}$.

\subsection{Proof of Theorem~\ref{T: main 2}, part (ii)}\label{SS: relationship}{\quad}\vspace{6pt}

By our usual arguments about M\"{o}bius-invariance, it suffices to deal with the case of subsets $E \subseteq F \subseteq \mathbb{R}^m$. We have to show that the following are equivalent:
\begin{enumerate}
\item $E \in \textup{CL}(m)$, $F$ is closed, and 
\[ 
F^\circ \subseteq \overline{E} \quad \text{and} \quad E \subseteq F\,;
\]
\item there exists a standard sequence $\left(w_n\right)$ in $\mathbb{H}^{m+1}$ such that 
\[ 
F = \Lambda\left(w_n\right) \quad \textup{ and } \quad E = \Lambda_c\left(w_n\right)\,.
\]
\end{enumerate} 
\begin{proof}
We saw in Lemma~\ref{L: dense ordinary implies uncountable conical} that $\Lambda\left(w_n\right)^\circ \subseteq \overline{\Lambda_c\left(w_n\right)}$, and the rest of the implication $\textup{(ii)} \Rightarrow \textup{(i)}$ is straightforward.

To prove $\textup{(i)} \Rightarrow \textup{(ii)}$ we use the notion of the dune $\mathcal{D}(U)$ of an open set $U \subseteq \mathbb{R}^m$, that was defined in the proof of Lemma~\ref{L: intersectOpens}.
  It is easy to check that the boundary of $\mathcal{D}\left(\mathbb{R}^m \setminus F\right)$ in $\mathbb{H}^{m+1}$ has  limit set equal to $\partial F$ but has empty conical limit set. Let $W$ be a subset of $\mathbb{H}^{m+1}$ such that $E=\Lambda_c(W)$. We define 
  \[ W' \, = \,  \left(W \setminus \mathcal{D}\left(\mathbb{R}^m \setminus F\right) \right)\, \cup \, \partial \mathcal{D}\left(\mathbb{R}^m \setminus F\right)\,.\]
 Then \[\Lambda(W')\, =\, \left(\,\overline{E} \setminus \left(\mathbb{R}^m \setminus F\right)\right) \cup \partial F \,=\, F\,,\] and
 \[ \Lambda_c(W') \,  = \, \left( E \setminus\left(\mathbb{R}^m \setminus F \right) \right) \cup \emptyset \, = \, E \,.\] 
Finally, as in the proof of Lemma~\ref{L: escaping and standard sequences} we can pick a standard sequence $(w_n)$ with the same conical limit set and  limit set as $W'$.
\end{proof}

\subsection{Proof of Theorem~\ref{T: main 2}, part (iii)}\label{SS: quasi}{\quad}\vspace{6pt}

In this section we prove that the class $\textup{CL}(m)$ is invariant under quasiconformal homeomorphisms of $S^m$, for $m \geq 2$.  The analogous result for $m=1$ is that that the class $\textup{CL}(1)$ is invariant under quasisymmetric homeomorphisms of $S^1$; the proof is similar and we omit it.  For the necessary background on quasisymmetric and quasiconformal mappings, see \cite{AV} or \cite{Vu}. Because conical limit sets are locally determined, and locally quasiconformal homeomorphisms are locally extensible to globally quasiconformal homeormorphisms, we can deduce from Theorem~\ref{T: main 2} (iii)  that $\textup{CL}(m)$ is invariant under locally quasiconformal homeomorphisms. 

It is possible to prove Theorem~\ref{T: main 2} (iii) directly from the definition of a locally quasisymmetric mapping. However, our theme is to think about conical limit sets as being subsets of the ideal boundary of hyperbolic space, and we therefore choose to present instead a proof that uses two well-known results about quasiconformal  homeomorphisms of $\partial \mathbb{B}^m$ and quasi-isometric homeomorphisms of $\mathbb{B}^m$ with the hyperbolic metric.

\begin{theorem*}[Tukia and V\"ais\"al\"a, \cite{TV2}]{\quad}\\
For $m \geq 2$, any quasiconformal homeomorphism $f: \mathbb{S}^{m} \to \mathbb{S}^{m}$ extends continuously to a quasiconformal homeomorphism $\hat{f}: \overline{\mathbb{B}^{m+1}} \to \overline{\mathbb{B}^{m+1}}$, that is bi-Lipschitz with respect to the hyperbolic metric on $\mathbb{B}^{m+1}$, with a constant that depends only on the constant of quasiconformality of the mapping $f$.
\end{theorem*}

\begin{theorem*}[Efremovic and Tihomirova, \cite{ET}] \label{T: quasi-geodesics}{\quad}\\
For $m\geq 1$, let $\hat{f}: \mathbb{B}^{m+1} \to \mathbb{B}^{m+1}$ be a quasi-isometry with respect to the hyperbolic metric, and let $x$ be a point in $\mathbb{S}^m$. The image curve $f([\jj,x])$ lies within a fixed distance of a geodesic ray, and $f([\jj,x])$ lands at a unique point of $\mathbb{S}^m$.
\end{theorem*}

Let $X \subseteq \mathbb{S}^m$ be a conical limit set and suppose that $f: \mathbb{S}^m \to \mathbb{S}^m$ is a quasiconformal homeomorphism. We have to show that the image $f(X)$ is also a conical limit set.  We have $X = \Lambda_c(w_n)$ for some escaping sequence of points $(w_n)$ in $\mathbb{B}^{m+1}$. Let $\hat{f}$ be a bi-Lipschitz extension of $f$ to $\overline{\mathbb{B}^{m+1}}$, guaranteed to exist by the above theorem of Tukia and V\"ais\"al\"a,  and let $k$ be the bi-Lipschitz constant of $\hat{f}$. Note that $\hat{f}$ is proper (since $\hat{f}^{-1}$ exists and is continuous), so that the sequence $(\hat{f}(w_{n}))$ is an escaping sequence. We now claim that if a subsequence $(w_{n_k})$ converges conically to $x \in X$ then $(\hat{f}(w_{n_k}))$ converges conically to $f(x) \in f(X)$.  Suppose that $(w_{n_k})$ converges conically to $x \in X$, so there is some $\alpha < \infty $ such that $\rho(w_{n_k}, [\jj,x]) \leq \alpha$ for all $k$. The image curve $\hat{f}([\jj,x])$ lands at $f(x)$ and remains within some distance $\delta>0$ of a geodesic ray $\gamma$.  We have
\[\rho(\hat{f}(w_{n_k}), \gamma) \,\leq\, k \rho(w_{n_k}, [\jj,x]) \, + \, \delta\, \le k\alpha + \delta.\]
Thus the sequence $(\hat{f}(w_{n_k}))$ converges conically to $f(x)$. Since $\hat{f}$ is a bi-Lipschitz homeomorphism extending $f$, we can run the same argument with $\hat{f}^{-1}$ to prove that $(w_{n_k})$ converges conically to $x$ only if $(\hat{f}(w_{n_k}))$ converges conically to $f(x)$. Hence $f(X) = \Lambda_c(\hat{f}(w_n))$.


\section{Open problems}\label{S: intersections}{\quad}\vspace{6pt}
In \cite{PT}, Piranian and Thron asked whether the intersection of two sets of divergence must also be a set of divergence, and this remains unsolved. We repeat the problem in general dimension:

\begin{problem} \label{P: closure under finite intersections}
Suppose $A, B \subseteq \mathbb{R}^m$ belong to $\textup{CL}(m)$. Must $A \cap B$ belong to $\textup{CL}(m)$?
\end{problem}

Here are two closely related problems:

\begin{problem} \label{P: intersections with G deltas}
Suppose $A, B \subseteq \mathbb{R}^m$, $A \in \textup{CL}(m)$, and $B$ is a $G_\delta$ set. Must $A \cap B$ belong to $\textup{CL}(m)$?
\end{problem}

\begin{problem} \label{P: closure under Cartesian products}
Suppose $X \subseteq \mathbb{R}^{m_1}$, $Y \subseteq \mathbb{R}^{m_2}$, $X \in \textup{CL}(m_1)$ and $Y \in \textup{CL}(m_2)$. Must  $X \times Y$ belong to $\textup{CL}(m_1 + m_2)$?
\end{problem}

The class $\text{CL}(m)$ is not closed under countable intersections, at least when $m>1$. To see this, choose any $G_{\delta\sigma\delta}$ subset $A$ of a codimension one sphere  $\mathbb{S}^{m-1} \subseteq\mathbb{S}^m$ such that $A$ is not a $G_{\delta\sigma}$ set. Express $A$ as $\bigcap_{n=1}^\infty A_n$, where the $A_n$ are $G_{\delta\sigma}$ subsets of $\mathbb{S}^{m-1}$. Theorem~\ref{T: main 1}(ii) shows that $A\notin \text{CL}(m)$. On the other hand, Theorem~\ref{T: main 1}(iv) shows that $A_n\in \text{CL}(m)$, for each $n$.

Our guess for Problem~\ref{P: closure under finite intersections} is that $A \cap B$ need not necessarily belong to $\textup{CL}(m)$. For example, consider the following countable subset of $\mathbb{H}^2$, where $\alpha > 2$ is a parameter.
\[ 
J(\alpha) = \left\{\left(\frac{1}{q^\alpha}, \frac{a}{q}\right)\, : \,a, q \in \mathbb{N}, \, 0 < a < q , \, \gcd(a,q) = 1 \right\}. 
\]
Then $\theta \in \Lambda_c(J(\alpha))$ if and only if $\theta \in [0,1]$ and $\theta$ can be approximated by rationals to order $\alpha$, that is, for some constant $c>0$ there exist infinitely many distinct rational numbers $a/q$ such that $| \theta - a/q | \le  c/q^\alpha$. Let $T: x \mapsto x+ \pi -3$ and consider the set $\Lambda_c(J(3)) \cap T\left(\Lambda_c(J(3)\right)$. This set is uncountable and dense in the interval $[\pi - 3, 1]$. It seems unlikely that this set should belong to $\textup{CL}(1)$: we see no reason to expect that both $\theta$ and $\theta - \pi$ can be approximated by rationals to order $3$ if and only if $\theta$ can be approximated at a certain rate by real numbers of some particular form.

A positive answer to Problem~\ref{P: closure under finite intersections} implies a positive answer to Problem \ref{P: intersections with G deltas} since every $G_{\delta}$ set is a conical limit set. Likewise a positive answer to Problem~\ref{P: closure under finite intersections} implies a positive answer to Problem~\ref{P: closure under Cartesian products}. Indeed, suppose $X = \Lambda_c(E)$ where $E \subseteq \mathbb{H}^{m_1+1}$. Then $X \times \mathbb{R}^{m_2} = \Lambda_c(E \times \mathbb{R}^{m_2})$, where we identify $\mathbb{H}^{m_1+1} \times \mathbb{R}^{m_2}$ with $\mathbb{H}^{m_1 + m_2+1}$ in the obvious way. Likewise $\mathbb{R}^{m_1} \times Y \in \textup{CL}(m_1 + m_2)$,  so  $X \times Y$ is the intersection of two conical limit sets.

 All the explicit examples that the authors know of subsets of $\mathbb{R}^m$ which are $G_{\delta \sigma}$ but not conical limit sets have been constructed from countable sets that are not conical limit sets. We have seen that any subset of a countable conical limit set is again a conical limit set. So we cannot hope to find a counterexample for Problem~\ref{P: closure under finite intersections} in which either $A$ or $B$ is countable. We now show that we cannot even hope to find a counterexample in which $A \cap B$ is countable: just take $A_1 = A$, $A_2 = B$ and $A_3 = A_4 =  \dots = \mathbb{R}^m$ in the following lemma. 

\begin{lemma}\label{L: countable countable intersection of CLs is CL}
If $A_1, A_2, \dots \in \textup{CL}(m)$ and $\bigcap_{i=1}^\infty A_i$ is countable, then $\bigcap_{i=1}^\infty A_i \, \in \textup{CL}(m)$. 
\end{lemma}
\begin{proof}
 Suppose for a contradiction that $\bigcap_{i=1}^\infty A_i$ is countable but not a conical limit set. For each $i$ there is a sequence $(w_n^{(i)})$ in $\mathbb{H}^{m+1}$ such that $A_i = \Lambda_c(w_n^{(i)})$. From the characterisation of conical limit sets among countable sets, we know that $\bigcap_{i=1}^\infty A_i$ has a non-empty subset $E$ such that $E \subseteq \textup{gd}(E)$. Now, each $A_i$ contains $E$, so we can find an uncountable set contained in the 1-conical limit set of each sequence $(w_n^{(i)})$ by modifying the proof of Lemma~\ref{L: E is uncountable}. We do this by substituting the sequence $(w_n^{(i(k))})$, $n=1,2,\dotsc$, for the sequence $(w_n)$ from Lemma~\ref{L: E is uncountable} at the $k^{th}$ step of the recursion, where $i(k)$ is a sequence that takes every positive integer value infinitely often. So each point of the resulting uncountable limit set is the $1$-conical limit of a sequence which contains an infinite subsequence of each sequence $(w_n^{(i)})$. This means that $\bigcap_{i=1}^\infty A_i$ contains $E$ and is therefore uncountable, which is a contradiction.
\end{proof} 
 
 If $X$ is a countable dense subset of $\mathbb{R}$, then $X \times \mathbb{R}$ is a not a conical limit set in $\mathbb{R}^2$. However, Lemma~\ref{L: countable countable intersection of CLs is CL} shows that it is not the intersection of countably many conical limit sets in $\mathbb{R}^2$. If it were, then so would be $(X \times \mathbb{R}) \cap (\mathbb{R} \times X) = X \times X$, which is a dense countable set in $\mathbb{R}^2$ and therefore not a conical limit set, contrary to the lemma. More generally, any set whose intersection with countably many quasisymmetric images of itself is countable but dense cannot be a conical limit set. This gives us another way to prove that a set is not a conical limit set, but it cannot help us to find a counterexample for Problem~\ref{P: closure under finite intersections} because no such set can arise as the intersection of finitely many conical limit sets.


\end{document}